# FLEXIBLE COVARIANCE ESTIMATION IN GRAPHICAL GAUSSIAN MODELS

By Bala Rajaratnam, Hélène Massam[1] and Carlos M. Carvalho

*Stanford University, York University and University of Chicago*

In this paper, we propose a class of Bayes estimators for the covariance matrix of graphical Gaussian models Markov with respect to a decomposable graph $G$. Working with the $W_{P_G}$ family defined by Letac and Massam [*Ann. Statist.* **35** (2007) 1278–1323] we derive closed-form expressions for Bayes estimators under the entropy and squared-error losses. The $W_{P_G}$ family includes the classical inverse of the hyper inverse Wishart but has many more shape parameters, thus allowing for flexibility in differentially shrinking various parts of the covariance matrix. Moreover, using this family avoids recourse to MCMC, often infeasible in high-dimensional problems. We illustrate the performance of our estimators through a collection of numerical examples where we explore frequentist risk properties and the efficacy of graphs in the estimation of high-dimensional covariance structures.

**1. Introduction.** In this paper we consider the problem of estimation of the covariance matrix $\Sigma$ of an $r$-dimensional graphical Gaussian model. Since the work of Stein [35] the problem of estimating $\Sigma$ is recognized as highly challenging. In recent years, the availability of high-throughput data from genomic, finance, marketing (among others) applications has pushed this problem to an extreme where, in many situations, the number of samples ($n$) is often much smaller than the number of parameters. When $n < r$ the sample covariance matrix $S$ is not positive definite but even when $n > r$, the eigenstructure tends to be systematically distorted unless $r/n$ is extremely small (see [12, 35]). Numerous papers have explored better alternative estimators for $\Sigma$ (or $\Sigma^{-1}$) in both the frequentist and Bayesian frameworks (see [4, 8, 9, 15, 16, 17, 19, 25, 26, 29, 35, 37]). Many of these estimators give

---

Received June 2007; revised February 2008.

[1]Supported by NSERC Discovery Grant A8946.

substantial risk reductions compared to the sample covariance estimator $S$ in small sample sizes. A common underlying property of many of these estimators is that they are shrinkage estimators in the sense of James–Stein [19, 34]. In particular the Bayesian approach often yields estimators which "shrink" toward a structure associated with a prespecified prior. One of the first papers to exploit this idea is [4] which shows that if the prior used on $\Sigma^{-1}$ is the standard conjugate, that is, a Wishart distribution, then for an appropriate choice of the shape (or shrinkage) and scale hyperparameters, the posterior mean for $\Sigma$ is a linear combination of $S$ and the prior mean (see Section 3.1). It is easy to show [see (3.16)] that the eigenvalues of such estimators are also shrinkage estimators of the eigenvalues of $\Sigma$. More recently, for high-dimensional complex datasets with $r$ often larger than $n$, regularization methods have been proposed, which impose structure on the estimators through zeros in the covariance or the precision matrix (see [2, 18, 30]). The idea of imposing zeros in the precision matrix is not new, however, and was introduced in [12] in a pioneering paper on covariance selection models which are particular cases of graphical Gaussian models. Graphical Gaussian models have proven to be excellent tools for the analysis of complex high-dimensional data where dependencies between variables are expressed by means of a graph [3, 21].

In this paper we combine the regularization approach given by graphical models with the Bayesian approach of shrinking toward a structure. Through a decision-theoretic approach, we derive Bayes estimators of the covariance and precision matrices under certain priors and given loss functions, such that the precision matrix has a given pattern of zeros. Indeed, we work within the context of graphical Gaussian models Markov with respect to a decomposable graph $G$. Restricting ourselves to decomposable graphs allows us to use the family of inverse $W_{P_G}$ Wishart distributions [27] as priors for $\Sigma$. This is a family of conjugate prior distributions for $\Sigma^{-1}$ which includes the Wishart when $G$ is complete (i.e., when the model is saturated) and the inverse of the hyper inverse Wishart, the current standard conjugate prior for $\Sigma^{-1}$, when the model is Markov with respect to $G$ decomposable. A potentially restrictive feature of the inverse of the hyper inverse Wishart (and the Wishart) is the fact that it has only one shape parameter. The family of $W_{P_G}$ Wishart distributions considered here has three important characteristics. First, it has $k+1$ shape parameters where $k$ is the number of cliques in $G$. Second, it forms a conjugate family with an analytically explicit normalizing constant. Third, the Bayes estimators can be obtained in closed-form.

In Section 2, we give some fundamentals of graphical models. In Section 3, we recall the properties of the $W_{P_G}$ family and its inverse, the $IW_{P_G}$, and we derive the mathematical objects needed for our estimators, that is, the explicit expression for the mean of the $IW_{P_G}$. Parallel to the development of the $IW_{P_G}$, we present in Section 4 a noninformative reference prior for $\Sigma$ (and the precision matrix $\Omega$). While offering an objective procedure that avoids the specification of hyperparameters, the reference prior also allows for closed-form posterior estimation as the posterior for $\Sigma$ remarkably falls within the $IW_{P_G}$ family. In Section 5, we derive the Bayes estimator under two commonly used loss functions adapted to graphical models and the prior considered in Sections 3 and 4. Finally, in Sections 6 and 7 we compare the performance of our estimators in a series of high-dimensional examples.

**2. Preliminaries.** Let $G = (V, E)$ be an undirected graph with vertex set $V = \{1, \ldots, r\}$ and edge-set $E$. Vertices $i$ and $j$ are said to be neighbors in $G$ if $(i, j) \in E$. Henceforth in this paper, we will assume that $G$ is decomposable [24], where a perfect order of the cliques is available. For $(C_1, \ldots, C_k)$ in a perfect order, we use the notation $H_1 = R_1 = C_1$ while for $j = 2, \ldots, k$ we write

$$H_j = C_1 \cup \cdots \cup C_j, \qquad R_j = C_j \setminus H_{j-1}, \qquad S_j = H_{j-1} \cap C_j.$$

The $S_j, j = 2, \ldots, k$ are the minimal separators of $G$. Some of these separators can be identical. We let $k' \leq k - 1$ denote the number of distinct separators and $\nu(S)$ denote the multiplicity of $S$, that is, the number of $j$ such that $S_j = S$. Generally, we will denote by $\mathcal{C}$ the set of cliques of a graph $G$ and by $\mathcal{S}$ its set of separators.

An $r$-dimensional Gaussian model is said to be Markov with respect to $G$ if for any edge $(i, j)$ not in $E$, the $i$th and $j$th variables are conditionally independent given all the other variables. Such models are known as covariance selection models [12] or graphical Gaussian models (see [24, 36]). Without loss of generality, we can assume that these models have zero mean and are characterized by the parameter set $P_G$ of positive definite precision (or inverse covariance) matrices $\Omega$ such that $\Omega_{ij} = 0$ whenever the edge $(i, j)$ is not in $E$. Equivalently, if we denote by $M$ the linear space of symmetric matrices of order $r$, by $M_r^+ \subset M$ the cone of positive definite (abbreviated $> 0$) matrices, by $I_G$ the linear space of symmetric incomplete matrices $x$ with missing entries $x_{ij}, (i, j) \notin E$ and by $\kappa \colon M \mapsto I_G$ the projection of $M$ into $I_G$, the parameter set of the Gaussian model can be described as the set of incomplete matrices $\Sigma = \kappa(\Omega^{-1}), \Omega \in P_G$. Indeed it is easy to verify that the entries $\Sigma_{ij}, (i, j) \notin E$ are such that

$$(2.1) \qquad \Sigma_{ij} = \Sigma_{i, V \setminus \{i, j\}} \Sigma_{V \setminus \{i, j\}, V \setminus \{i, j\}}^{-1} \Sigma_{V \setminus \{i, j\}, j},$$

and are therefore not free parameters of the Gaussian models. We are therefore led to consider the two cones

$$(2.2) \qquad P_G = \{y \in M_r^+ \,|\, y_{ij} = 0, (i, j) \notin E\},$$

$$(2.3) \qquad Q_G = \{x \in I_G \,|\, x_{C_i} > 0, i = 1, \ldots, k\},$$

where $P_G \subset Z_G$ and $Q_G \subset I_G$, where $Z_G$ denotes the linear space of symmetric matrices with zero entries $y_{ij}, (i,j) \notin E$.

Gröne et al. [14] proved the following:

PROPOSITION 2.1. *When $G$ is decomposable, for any $x$ in $Q_G$ there exists a unique $\hat{x}$ in $M_r^+$ such that for all $(i,j)$ in $E$ we have $x_{ij} = \hat{x}_{ij}$ and such that $\hat{x}^{-1}$ is in $P_G$.*

This defines a bijection between $P_G$ and $Q_G$:

$$\varphi: y = (\hat{x})^{-1} \in P_G \mapsto x = \varphi(y) = \kappa(y^{-1}) \in Q_G, \tag{2.4}$$

where $\kappa$ denotes the projection of $M$ into $I_G$.

If for any complete subset $A \subseteq V$, $x_A = (x_{ij})_{i,j \in A}$ is a matrix and we denote by $(x_A)^0 = (x_{ij})_{i,j \in V}$ the matrix such that $x_{ij} = 0$ for $(i,j) \notin A \times A$, then the explicit expression of $\hat{x}^{-1}$ is

$$y = \hat{x}^{-1} = \sum_{C \in \mathcal{C}} (x_C^{-1})^0 - \sum_{S \in \mathcal{S}} \nu(S)(x_S^{-1})^0. \tag{2.5}$$

For $(x,y) \in I_G \times Z_G$, we define the notion of trace as follows:

$$\mathrm{tr}(xy) = \langle x, y \rangle = \sum_{(i,j) \in E} x_{ij} y_{ij}. \tag{2.6}$$

Note that for $x \in Q_G$ and $y \in P_G$, $\langle x, y \rangle = \mathrm{tr}(\hat{x}y)$, where $\mathrm{tr}(\hat{x}y)$ is defined in the classical way. In the sequel, we will also need the following. If, for $y \in P_G$ we write $y = \hat{\sigma}^{-1}$ with $\sigma \in Q_G$, we have, for $x \in Q_G$, the two formulas

$$\langle x, \hat{\sigma}^{-1} \rangle = \sum_{C \in \mathcal{C}} \langle x_C, \sigma_C^{-1} \rangle - \sum_{S \in \mathcal{S}} \nu(S) \langle x_S, \sigma_S^{-1} \rangle, \tag{2.7}$$

$$\det \hat{x} = \frac{\prod_{C \in \mathcal{C}} (\det x_C)}{\prod_{S \in \mathcal{S}} (\det x_S)^{\nu(S)}}. \tag{2.8}$$

The graphical Gaussian model Markov with respect to $G$ is therefore the family of distributions

$$\mathcal{N}_G = \{N_r(0, \Sigma), \Sigma \in Q_G\} = \{N_r(0, \Sigma), \Omega = \hat{\Sigma}^{-1} \in P_G\}.$$

In this paper, we will study various estimators of $\Sigma \in Q_G$ and $\Omega \in P_G$. We will write $mle$ and $mle_g$ for "maximum likelihood estimate" in the saturated model and in the graphical model, respectively. Also, in this paper, we will use the general symbol $\tilde{\theta}$ to denote an estimator of $\theta$ rather than the traditional $\hat{\theta}$ as the notation $\hat{\theta}$ has been reserved for the completion process (see Proposition 2.1). The $mle_g$, $\tilde{\Omega}_g$, for the parameter $\Omega \in P_G$ in $\mathcal{N}_G$ is well known (see [24], page 138). If $Z_i, i = 1, \ldots, n$ is a sample from the $N_r(0, \Sigma)$ distribution in $\mathcal{N}_G$, if we write $U = \sum_{i=1}^{n} Z_i Z_i^t$ and $S = \frac{U}{n}$ and if $n > \max_{C \in \mathcal{C}} |C|$, then $\widehat{\Omega}_g$ exists and is equal to

$$(2.9) \qquad \widehat{\Omega}_g = \sum_{C \in \mathcal{C}} (S_C^{-1})^0 - \sum_{S \in \mathcal{S}} \nu(S)(S_S^{-1})^0,$$

where clearly $S$ as a subscript or $S$ in $\nu(S)$ refers to the separator while the remaining $S$'s refer to the sample covariance matrix. If we assume that the graph is saturated, then, clearly the *mle* is $\tilde{\Omega} = S^{-1}$.

Finally, we need to recall some standard notation for various block submatrices: for $x \in Q_G$, $x_{C_j}, j = 1, \ldots, k$ are well defined and for $j = 2, \ldots, k$, it will be convenient to use the following:

$$
(2.10) \qquad
\begin{aligned}
x_{S_j} &= x_{\langle j \rangle}, & x_{R_j, S_j} &= x_{[j \rangle} = x_{\langle j]}^t, \\
x_{[j]} &= x_{R_j}, & x_{[j] \cdot} &= x_{[j]} - x_{[j \rangle} x_{\langle j \rangle}^{-1} x_{\langle j]},
\end{aligned}
$$

where $x_{\langle j \rangle} \in M_{s_j}^+, x_{[j]} \in M_{c_j - s_j}^+, x_{[j \rangle} \in L(\mathbb{R}^{c_j - s_j}, \mathbb{R}^{s_j})$, the set of linear applications from $\mathbb{R}^{c_j - s_j}$ to $\mathbb{R}^{s_j}$. We will also use the notation $x_{[12 \rangle}$ and $x_{[1] \cdot}$ for

$$x_{[12 \rangle} = x_{C_1 \setminus S_2, S_2} x_{S_2}^{-1} \quad \text{and} \quad x_{[1] \cdot} = x_{C_1 \setminus S_2 \cdot S_2} = x_{C_1 \setminus S_2} - x_{C_1 \setminus S_2, S_2} x_{S_2}^{-1} x_{S_2, C_1 \setminus S_2}.$$

## 3. Flexible conjugate priors for $\Sigma$ and $\Omega$.

When the Gaussian model is saturated, that is, $G$ is complete, the conjugate prior for $\Omega$, as defined by Diaconis and Ylvisaker [13] (henceforth abbreviated DY) is the Wishart distribution. The induced prior for $\Sigma$ is then the inverse Wishart $IW_r(p, \theta)$ with density

$$(3.1) \qquad IW_r(p, \theta; dx) = \frac{|\theta|^p}{\Gamma_r(p)} |x|^{-p - (r+1)/2} \exp -\langle \theta, x^{-1} \rangle \mathbf{1}_{M^+}(x) \, dx,$$

where $p > \frac{r-1}{2}$ is the shape parameter, $\Gamma_r(p)$ the multivariate gamma function (as given on page 61 of [31]) and $\theta \in M^+$ is the scale parameter.

As we have seen in the previous section, when $G$ is not complete, $M^+$ is no longer the parameter set for $\Sigma$ or the parameter set for $\Omega$. The DY conjugate prior for $\Sigma \in Q_G$ was derived by [11] and called the hyper inverse Wishart ($HIW$). The induced prior for $\Omega \in P_G$ was derived by [32] and we will call it the $G$-Wishart. The $G$-Wishart and the hyper inverse Wishart are certainly defined on the right cones but they essentially have the same type of parametrization as the Wishart with a scale parameter $\theta \in Q_G$ and a one-dimensional shape parameter $\delta$.

3.1. *The $W_{P_G}$ distribution and its inverse.* Letac and Massam [27] introduced a new family of conjugate priors for $\Omega \in P_G$ with a $k+1$-dimensional shape parameter thus leading to a richer family of priors for $\Omega$, and therefore for $\Sigma = \kappa(\Omega^{-1}) \in Q_G$ through the induced prior. It is called the type II Wishart family. Here, we prefer to call it the family of $W_{P_G}$-Wishart distributions in order to emphasize that it is defined on $P_G$. Details of this distribution can be found in Section 3 of [27]. We will first recall here some of its main features and then derive some new properties we shall need later in this paper. Let $\alpha$ and $\beta$ be two real-valued functions on the collection $\mathcal{C}$ and $\mathcal{S}$ of cliques and separators, respectively, such that $\alpha(C_i) = \alpha_i, \beta(S_j) = \beta_j$ with $\beta_i = \beta_j$ if $S_i = S_j$. Let $c_i = |C_i|$ and $s_i = |S_i|$ denote the cardinality of $C_i$ and $S_i$, respectively. The family of $W_{P_G}$-Wishart distributions is the natural exponential family generated by the measure $H_G(\alpha, \beta, \varphi(y))\nu_G(dy)$ on $P_G$ where $\varphi(y)$ is as defined in (2.4) and where, for $x \in Q_G$,

$$(3.2) \qquad H_G(\alpha, \beta; x) = \frac{\prod_{C \in \mathcal{C}}(\det x_C)^{\alpha(C)}}{\prod_{S \in \mathcal{S}}(\det x_S)^{\nu(S)\beta(S)}},$$

$$(3.3) \qquad \nu_G(dy) = H_G(\tfrac{1}{2}(c+1), \tfrac{1}{2}(s+1); \varphi(y)) \mathbf{1}_{P_G}(y)\, dy.$$

The parameters $(\alpha, \beta)$ are in the set $\mathcal{B}$ such that the normalizing constant is finite for all $\theta \in Q_G$ and such that it factorizes into the product of $H_G(\alpha, \beta; \theta)$ and a function $\Gamma_{II}(\alpha, \beta)$ of $(\alpha, \beta)$ only, given below in (3.5).

The set $\mathcal{B}$ is not known completely but we know that $\mathcal{B} \supseteq \bigcup_P B_P$ where, if, for each perfect order $P$ of the cliques of $G$, we write $J(P, S) = \{j = 2, \ldots, k | S_j = S\}$, then $B_P$ is the set of $(\alpha, \beta)$ such that:

1. $\sum_{j \in J(P,S)}(\alpha_j + \frac{1}{2}(c_j - s_j)) - \nu(S)\beta(S) = 0$, for all $S$ different from $S_2$;
2. $-\alpha_q - \frac{1}{2}(c_q - s_q - 1) > 0$ for all $q = 2, \ldots, k$ and $-\alpha_1 - \frac{1}{2}(c_1 - s_2 - 1) > 0$;
3. $-\alpha_1 - \frac{1}{2}(c_1 - s_2 + 1) - \gamma_2 > \frac{s_2 - 1}{2}$, where $\gamma_2 = \sum_{j \in J(P,S_2)}(\alpha_j - \beta_2 + \frac{c_j - s_2}{2})$.

As can be seen from the conditions above, the parameters $\beta(S), S \in \mathcal{S}$ are linked to the $\alpha(C), C \in \mathcal{C}$ by $k' - 1$ linear equalities and various linear inequalities and therefore $\mathcal{B}$ contains the set $B_P$ of dimension at least $k + 1$, for each perfect order $P$. We can now give the formal definition of the $W_{P_G}$ family.

DEFINITION 3.1. For $(\alpha, \beta) \in \mathcal{B}$, the $W_{P_G}$-Wishart family of distributions is the family $\mathcal{F}_{(\alpha,\beta), P_G} = \{W_{P_G}(\alpha, \beta, \theta; dy), \theta \in Q_G\}$ where

$$(3.4) \qquad W_{P_G}(\alpha, \beta, \theta; dy) = e^{-\langle \theta, y \rangle} \frac{H_G(\alpha, \beta; \varphi(y))}{\Gamma_{II}(\alpha, \beta)H_G(\alpha, \beta; \theta)} \nu_G(dy)$$

and

$$\Gamma_{II}(\alpha, \beta) = \pi^{((c_1 - s_2)s_2 + \sum_{j=2}^{k}(c_j - s_j)s_j)/2}$$

(3.5)
$$\times \Gamma_{s_2}\left[-\alpha_1 - \frac{c_1 - s_2}{2} - \gamma_2\right]\Gamma_{c_1 - s_2}(-\alpha_1)\prod_{j=2}^{k}\Gamma_{c_j - s_j}(-\alpha_j).$$

We can also, of course, define the inverse $W_{P_G}(\alpha, \beta, \theta)$ distribution as follows. If $Y \sim W_{P_G}(\alpha, \beta, \theta)$, then $X = \varphi(Y) \sim IW_{P_G}(\alpha, \beta, \theta)$ with distribution on $Q_G$ given by (see (3.8) in [27])

(3.6)
$$IW_{P_G}(\alpha, \beta, \theta; dx) = \frac{e^{-\langle\theta, \hat{x}^{-1}\rangle}H_G(\alpha, \beta; x)}{\Gamma_{II}(\alpha, \beta)H_G(\alpha, \beta; \theta)}\mu_G(dx),$$

(3.7)
$$\text{where } \mu_G(dx) = H_G(-\tfrac{1}{2}(c+1), -\tfrac{1}{2}(s+1); x)\mathbf{1}_{Q_G}(x)\, dx.$$

The hyper inverse Wishart is a special case of the $IW_{P_G}$ distribution for

(3.8)
$$\alpha_i = -\frac{\delta + c_i - 1}{2}, \qquad i = 1, \ldots, k,$$

$$\beta_i = -\frac{\delta + s_i - 1}{2}, \qquad i = 2, \ldots, k,$$

which are all functions of the same one-dimensional parameter $\delta$. It is traditional to denote the hyper inverse Wishart, that is, this particular $IW_{P_G}$, as the $HIW(\delta, \theta)$ and this is the notation we will use in Section 6.

Corollary 4.1 of [27] states that the $IW_{P_G}$ is a family of conjugate distributions for the scale parameter $\Sigma$ in $\mathcal{N}_G$; more precisely, we have:

PROPOSITION 3.1. *Let $G$ be decomposable and let $P$ be a perfect order of its cliques. Let $(Z_1, \ldots, Z_n)$ be a sample from the $N_r(0, \Sigma)$ distribution with $\Sigma \in Q_G$. If the prior distribution on $2\Sigma$ is $IW_{P_G}(\alpha, \beta, \theta)$ with $(\alpha, \beta) \in B_P$ and $\theta \in Q_G$, the posterior distribution of $2\Sigma$, given $nS = \sum_{i=1}^{n} Z_i Z_i^t$, is $IW_{P_G}(\alpha - \frac{n}{2}, \beta - \frac{n}{2}, \theta + \kappa(nS))$, where $\alpha - \frac{n}{2} = (\alpha_1 - \frac{n}{2}, \ldots, \alpha_k - \frac{n}{2})$ and $\beta - \frac{n}{2} = (\beta_2 - \frac{n}{2}, \ldots, \beta_k - \frac{n}{2})$ are such that $(\alpha - \frac{n}{2}, \beta - \frac{n}{2}) \in B_P$ and $\theta + \kappa(nS) \in Q_G$ so that the posterior distribution is well defined. Equivalently, we may say that if the prior distribution on $\frac{1}{2}\Omega$ is $W_{P_G}(\alpha, \beta, \theta)$, then the posterior distribution of $\frac{1}{2}\Omega$ is $W_{P_G}(\alpha - \frac{n}{2}, \beta - \frac{n}{2}, \theta + \kappa(nS))$.*

For the expression of the Bayes estimators we will give in Section 5, we need to know the explicit expression of the posterior mean of $\Omega$ and $\Sigma$ when the prior on $\Omega$ is the $W_{P_G}$ or equivalently when the prior on $\Sigma$ is the $IW_{P_G}$. The mean of the $W_{P_G}$ can be immediately obtained by differentiation of the cumulant generating function since the $W_{P_G}$ family is a natural exponential family. From (3.4), from Corollary 3.1 and from (4.25) in [27], we easily obtain the posterior mean for $\Omega = \widehat{\Sigma}^{-1}$ as follows.

PROPOSITION 3.2. *Let $S$ and $\Omega$ be as in Corollary 3.1; then the posterior mean of $\Omega$, given $nS$, is*

$$
\begin{aligned}
(3.9) \quad E(\Omega|S) = -2\Bigg[ &\sum_{j=1}^{k}\left(\alpha_j - \frac{n}{2}\right)((\theta + \kappa(nS))^{-1}_{C_j})^0 \\
&- \sum_{j=2}^{k}\left(\beta_j - \frac{n}{2}\right)((\theta + \kappa(nS))^{-1}_{S_j})^0 \Bigg].
\end{aligned}
$$

Since the $IW_{P_G}$ is not an exponential family, its expected value is not as straightforward to derive. It is given in the following theorem.

THEOREM 3.1. *Let $X$ be a random variable on $Q_G$ such that $X \sim IW_{P_G}(\alpha, \beta, \theta)$ with $(\alpha, \beta) \in B_P$ and $\theta \in Q_G$; then $E(X)$ is given by (3.10)–(3.14):*

$$
\begin{aligned}
(3.10) \quad E(x_{\langle 2 \rangle}) &= \frac{\theta_{\langle 2 \rangle}}{-(\alpha_1 + ((c_1 - s_2)/2) + \gamma_2) - ((s_2 + 1)/2)} \\
&= \frac{\theta_{\langle 2 \rangle}}{-(\alpha_1 + ((c_1 + 1)/2) + \gamma_2)},
\end{aligned}
$$

$$
(3.11) \quad E(x_{C_1 \setminus S_2, S_2}) = \frac{\theta_{C_1 \setminus S_2, S_2}}{-(\alpha_1 + ((c_1 + 1)/2) + \gamma_2)},
$$

$$
\begin{aligned}
(3.12) \quad E(x_{C_1 \setminus S_2}) &= \frac{\theta_{[1].}}{-(\alpha_1 + ((c_1 - s_2 + 1)/2))} \\
&\quad \times \left(1 - \frac{s_2}{2(\alpha_1 + ((c_1 + 1)/2) + \gamma_2)}\right) \\
&\quad + \frac{\theta_{C_1 \setminus S_2, S_2}\theta_{\langle 2 \rangle}^{-1}\theta_{S_2, C_1 \setminus S_2}}{-(\alpha_1 + ((c_1 + 1)/2) + \gamma_2)},
\end{aligned}
$$

*and for $j = 2, \ldots, k$*

$$
(3.13) \quad E(x_{[j]}) = E((x_{[j]}x_{\langle j \rangle}^{-1}))E(x_{\langle j \rangle}) = \theta_{[j]}\theta_{\langle j \rangle}^{-1}E(x_{\langle j \rangle}),
$$

$$
\begin{aligned}
(3.14) \quad E(x_{[j]}) &= \frac{\theta_{[j].}}{-(\alpha_j + ((c_j - s_j + 1)/2))}\left(1 + \frac{1}{2}\operatorname{tr}(\theta_{\langle j \rangle}^{-1}E(x_{\langle j \rangle}))\right) \\
&\quad + \theta_{[j]}\theta_{\langle j \rangle}^{-1}E(x_{\langle j \rangle})\theta_{\langle j \rangle}^{-1}\theta_{\langle j].}.
\end{aligned}
$$

The proof is rather long and technical and given in the Appendix. Let us note here that (3.10)–(3.14) can also be written in a closed-form expression of the Choleski type, that is, $E(X) = T^t D T$ with $T$ lower triangular and $D$ diagonal, where the shape parameters $(\alpha, \beta)$ are solely contained in $D$. We do not give it here for the sake of brevity.

The important consequence of this theorem is that, from (3.10) to (3.14), we can rebuild $E(X) \in Q_G$. Indeed, by definition of $Q_G$, $E(X)$ is made up first of $E(X_{C_1})$ which is given by (3.10), (3.11), and its transpose (3.12) and then, successively, of the $j$th "layer": $E(X_{[j]})$ and its transpose, and $E(X_{[j]})$, for each $j = 2, \ldots, k$. These are immediately obtained from (3.13) and (3.14) since, by definition, $S_j \subseteq H_{j-1}$ and therefore the quantity $E(X_{(j)})$ is a sub-block of $E(X_{H_{j-1}})$ and has therefore already been obtained in the first $j-1$ steps. We can therefore now deduce the posterior mean of $\Sigma$ when the prior is $IW_{P_G}(\alpha, \beta, \theta)$.

COROLLARY 3.1. *Let $S$ and $\Sigma$ be as in Corollary 3.1; then the posterior mean for $\Sigma$ when the prior distribution on $2\Sigma$ is $IW_{P_G}(\alpha, \beta, \theta)$ is given by (3.10)–(3.14) where $X$ is replaced by $2\Sigma$, $\theta$ is replaced by $\theta + \kappa(nS)$ and $(\alpha_i, \beta_i)$'s are replaced by $\alpha_i - \frac{n}{2}, \beta_i - \frac{n}{2}$'s.*

3.2. *Shrinkage by layers and the choice of the scale parameter $\theta$.* When we use the $IW_{P_G}(\alpha, \beta, \theta)$ as a prior distribution for the scale parameter $\Sigma$, we have to make a choice for the shape hyperparameters $(\alpha, \beta)$ and the scale hyperparameter $\theta$. When $G$ is complete, the $IW_{P_G}(\alpha, \beta, \theta)$ becomes the regular inverse Wishart $IW(p, \theta)$ as given in (3.1). When $G$ is decomposable and one uses the hyper inverse Wishart $HIW(\delta, \theta)$, in the absence of prior information, it is traditional to take $\theta$ to be equal to the identity or a multiple of the identity and $\delta$ small, such as 3, for example (see [21]).

The scale parameter, however, can play an important role if we have some prior knowledge on the structure of the covariance matrix (see [4]) and we are interested in "shrinking" the posterior mean of $\Sigma$ toward a given target.

In the saturated case, for a sample of size $n$ from the $N(0, \Sigma)$ distribution with a Wishart $W(\frac{\nu}{2}, (\nu D)^{-1})$ prior on $\Omega = \Sigma^{-1}$, the posterior mean of $\Sigma$ is

$$(3.15) \qquad E(\Sigma | S) = \frac{\nu D + nS}{\nu + n - r - 1}.$$

First, we note that when $n$ is held fixed and $\nu$ is allowed to grow, the posterior mean tends toward $D$ while if $\nu$ is held fixed and $n$ is allowed to grow, the estimator tends toward $S$. Next, let us consider the eigenvalues of the posterior mean. If we take $D = \bar{l} I$ where $\bar{l}$ is the average of the eigenvalues $l_1, \ldots, l_r$ of the *mle* $S$, then it is easy to see that the eigenvalues $g_i$, $i = 1, \ldots, r$ of $E(\Sigma | S)$ are

$$(3.16) \qquad g_i = \frac{\nu \bar{l} + n l_i}{\nu - (r+1) + n},$$

nearly a weighted average of $\bar{l}$ and $l_i$. Some simple algebra will show that for $l_i < \bar{l}$ we always have $l_i < g_i$ and that for $i$ such that $l_i > \bar{l}$, that is, for $C_i = \frac{l_i}{\bar{l}} > 1$, we will have $g_i < l_i$ whenever $\nu > \frac{C_i}{C_i - 1}(r + 1)$. Since in order for the prior to be proper, we must have that $\nu > r - 1$, we see that this condition is very weak as long as $\frac{C_i}{C_i - 1}$ is close to 1. When the condition $\nu > r - 1$ is satisfied, the eigenvalues of the posterior mean are shrunk toward $\bar{l}$ and the span of the eigenvalues of $E(\Sigma | S)$ is smaller than the span of the eigenvalues of $S$, which generally can be used to correct the instability of $S$. (We note that if $C_i = \frac{l_i}{\bar{l}}$ is sufficiently large, $\frac{C_i}{C_i - 1}$ will be sufficiently close to 1 and if $\frac{l_i}{\bar{l}}$ is close to 1, then there is really no need to shrink the eigenvalues.)

In Section 5 we show that our Bayes estimators can be expressed in terms of the posterior mean of $\Sigma$ and $\Omega$ with the $IW_{P_G}$ and the $W_{P_G}$, respectively, as priors. One would like to be able to prove properties for the eigenvalues of our estimators similar to those of the posterior mean under the Wishart in the saturated case. This is beyond the scope of this paper. However, we observe in the numerical examples given in Sections 6 and 7 that the eigenvalues of our estimators do have shrinkage properties. With this motivation, in Sections 6 and 7, we will use $IW_{P_G}$ priors with $\theta$ so that the prior mean of $\Sigma$ is the identity as well as with $\theta$ equal to the identity. Thus, we first derive $\theta$ so that $E(\Sigma) = \frac{1}{2}E(IW_{P_G}(\alpha, \beta, \theta)) = I$ and then, we will argue that our estimators can be viewed as shrinkage estimators in the sense of shrinkage toward structure.

LEMMA 3.1. *Let* $\Sigma \in Q_G$ *be such that* $2\Sigma \sim IW_{P_G}(\alpha, \beta, \theta)$ *for given* $(\alpha, \beta) \in \mathcal{A}$. *In order to have* $E(\Sigma) = I$ *it is sufficient to choose* $\theta$ *as a diagonal matrix with diagonal elements equal to*

$$\theta_{ll} = -2\left(\alpha_1 + \frac{c_1 - s_2 + 1}{2}\right)\left(1 - \frac{s_2}{2(\alpha_1 + ((c_1 + 1)/2) + \gamma_2)}\right)^{-1} \quad \text{for } l \in [1],$$

$$\theta_{ll} = -2\left(\alpha_1 + \frac{c_1 - s_2}{2} + \gamma_2\right) - (s_2 + 1) \quad \text{for } l \in \langle 2 \rangle,$$

$$\theta_{ll} = -2\left(\alpha_j + \frac{c_j - s_j + 1}{2}\right)\left(1 + \frac{1}{2}\operatorname{tr}(\theta_{\langle j \rangle}^{-1} E(x_{\langle j \rangle}))\right)^{-1}$$
$$\text{for } l \in [j], j = 2, \ldots, k.$$

The proof is immediate from (3.10)–(3.14).

Let us now argue that one of our estimators (to be derived in Section 5), $\widetilde{\Omega}_{L_1}^{W_{P_G}}$, equal to the inverse of the completion of the posterior mean $E(\Sigma | S)$ of $\Sigma$ when $\Sigma \sim IW_{P_G}(\alpha, \beta, \theta)$, can be viewed as a shrinkage estimator. It follows from Theorem 4.4 of [27] that, when the prior on $\Sigma$ is the $IW_{P_G}(\alpha, \beta, \theta)$, $\Sigma_{[i]} \sim IW_{c_i - s_i}(-\alpha_i, \theta_{[i]})$ as defined in (3.1). Then, since $nS_{H_i} \sim W_{|H_i|}(\frac{n}{2}, \Sigma_{H_i})$ and thus $nS_{[i]\cdot} \sim W_{c_i-s_i}(\frac{n-s_i}{2}, \Sigma_{[i]\cdot})$, through an argument parallel to the one for (3.15), it follows that the posterior mean $E(\Sigma_{[i]\cdot}|S)$ is a linear combination of $S_{[i]\cdot}$ and $\theta_{[i]\cdot}$, with $\alpha_i$ playing a role parallel to that of $\nu$ in (3.15), and is therefore a shrinkage estimator of $\Sigma_{[i]\cdot}$. Thus, we can shrink with different intensities various parts of the matrix $S$. The posterior mean $E(\Sigma|S)$ is reconstructed, layer by layer, as can be seen in the proof of Theorem 3.1 using $E(\Sigma_{[i]\cdot}|S)$ as a building block [see (A.4)] through what we might call a conditional Choleski reconstruction. The resulting estimate $\widetilde{\Omega}_{L_1}^{W_{P_G}}$ can therefore be regarded as a shrinkage estimator. A similar argument can be made for all our Bayes estimators.

**4. Reference prior.** In this section, we derive a reference prior for $\Sigma$ and therefore $\Omega$ (see [1, 5, 6]). This is done first by reparametrizing the density of the *mle* $\kappa(S)$ of $\Sigma$ with a parametrization naturally induced by a given perfect order $P$ of cliques. The parameters are, in fact, the elements of the Choleski decomposition of $\widehat{\Sigma}$ for the order of the vertices given by $P$. As we will see below, the density of the *mle* of $\Sigma$ belongs to a natural exponential family and we will therefore follow the method given by Datta and Ghosh [10] and later used by Consonni and Veronese [7] in the context of general Wishart distributions, to derive the reference prior for the new parameter. We will then consider the induced prior on $\Sigma$. It was shown in [27], equations (4.14)–(4.18), that if, for $X = \kappa(nS)$ and a given perfect order of the cliques $P$, we make the following change of variable:

$$x \mapsto \xi = (x_{[1]\cdot}, x_{[12]}, x_{\langle 2 \rangle}, x_{[j]\cdot}, x_{[j]}x_{\langle j \rangle}^{-1}, j = 2, \ldots, k),$$

then the density of the new variable $\Xi$ is

$$
\begin{aligned}
W_{Q_G}^{**}&(\alpha, \beta, \sigma; d\xi) \\
&\propto |\sigma_{[1]\cdot}^{-1}|^{p-(s_2/2)} |x_{[1]\cdot}|^{p-(s_2/2)-((c_1-s_2+1)/2)} e^{-\langle x_{[1]\cdot}, \sigma_{[1]\cdot}^{-1} \rangle} \, dx_{[1]\cdot} \\
&\quad \times |\sigma_{\langle 2 \rangle}^{-1}|^{p} |x_{\langle 2 \rangle}|^{p-((s_2+1)/2)} e^{-\langle x_{\langle 2 \rangle}, \sigma_{S_2}^{-1} \rangle} \\
&\quad \times |\sigma_{[1]\cdot}^{-1}|^{+(s_2/2)} |x_{\langle 2 \rangle}|^{(c_1-s_2)/2} \\
&\quad \times (e^{-\langle (x_{[12]}-\sigma_{[12]}), \sigma_{[1]\cdot}^{-1}(x_{[1,2]}-\sigma_{[1,2]})x_{\langle 2 \rangle} \rangle} \, dx_{[12]}) \, dx_{\langle 2 \rangle} \\
&\quad \times \prod_{j=2}^{k} |\sigma_{[j]\cdot}^{-1}|^{+(s_j/2)} |x_{\langle j \rangle}|^{(c_j-s_j)/2} \\
&\qquad \times e^{-\langle (x_{[j]}x_{\langle j \rangle}^{-1}-\sigma_{[j]}\sigma_{\langle j \rangle}^{-1}), \sigma_{[j]\cdot}^{-1}(x_{[j]}x_{\langle j \rangle}^{-1}-\sigma_{[j]}\sigma_{\langle j \rangle}^{-1})x_{\langle j \rangle} \rangle}
\end{aligned}
\tag{4.1}
$$

$$\times \, |\sigma_{[j] \cdot}^{-1}|^{p-(s_j/2)} |x_{[j] \cdot}|^{p-(s_j/2)-((c_j-s_j+1)/2)}$$

$$\times \, e^{-\langle x_{[j] \cdot}, \sigma_{[j] \cdot}^{-1} \rangle} \, dx_{C_1} \prod_{j=2}^{k} d(x_{[j]} x_{\langle j \rangle}^{-1}) \, dx_{[j] \cdot},$$

and the new parametrization replacing $\sigma$, in the order induced by the order of the new variables, is clearly

$$(4.2) \qquad \phi = (\sigma_{\langle 2 \rangle}, (\sigma_{[1] \cdot}^{-1}, \sigma_{[12]}), (\sigma_{[j] \cdot}^{-1}, \sigma_{[j]} \sigma_{\langle j \rangle}^{-1}), j = 2, \dots, k).$$

We now derive the reference prior for $\phi$ and the induced prior for $\sigma$.

THEOREM 4.1.   *Consider the scale parameter $\Sigma \in Q_G$ for the Gaussian model $\mathcal{N}_G$. Let $\sigma = 2\Sigma$ and let $\phi$ be the ordered parameter as defined in (4.2). The reference prior for $\phi$ is independent of the order of the components and has density equal to*

$$(4.3) \qquad \pi^{\phi}(\phi) = |\sigma_{[1] \cdot}|^{((c_1+1)/2)-s_2} |\sigma_{\langle 2 \rangle}|^{(s_2+1)/2} \prod_{j=2}^{k} |\sigma_{[j] \cdot}|^{((c_j+1)/2)-s_j}.$$

*Moreover the induced reference prior for the parameter $\sigma \in Q_G$ is*

$$(4.4) \qquad \pi^{\sigma}(\sigma) = \frac{|\sigma_{C_1}|^{-(c_1+1)/2} \prod_{j=2}^{k} |\sigma_{C_j}|^{-(c_j+1)/2}}{|\sigma_{S_2}|^{((c_1+c_2)/2)-s_2-((s_2+1)/2)} \prod_{j=3}^{k} |\sigma_{S_j}|^{((c_j-s_j)/2)-((s_j+1)/2)}},$$

*which corresponds to an improper $IW_{P_G}(\alpha, \beta, 0)$ distribution with*

$$(4.5) \qquad \begin{aligned} &\alpha_j = 0, \qquad j = 1, \dots, k, \\ &\beta_2 = \frac{c_1 + c_2}{2} - s_2, \qquad \beta_j = \frac{c_j - s_j}{2}, \qquad j = 3, \dots, k. \end{aligned}$$

The proof of the theorem is given in the Appendix. Let us note here that the fact that the induced prior on $\sigma$ is an $IW_{P_G}$, albeit an improper one, is not too surprising since (see Theorem 4.1 of [27]) the $IW_{P_G}$ is a conjugate distribution for the scale parameter of the distribution of $\kappa(nS)$.

Because the distribution (4.4) has the form of an $IW_{P_G}$, its posterior given $U = nS$ is an $IW_{P_G}$ with parameters

$$(4.6) \qquad \begin{aligned} &\alpha_j = 0 - \frac{n}{2}, \qquad j = 1, \dots, k, \\ &\beta_2 = \frac{c_1 + c_2}{2} - s_2 - \frac{n}{2}, \qquad \beta_j = \frac{c_j - s_j}{2} - \frac{n}{2}, \qquad j = 3, \dots, k \end{aligned}$$

and

$$\theta = \kappa(nS).$$

Of course, $\kappa(ns) \in Q_G$ and it is easy to check that $(\alpha, \beta) \in \mathcal{B}$, that is, the posterior distribution is a proper $IW_{P_G}$. As in Section 3, we now need to compute the explicit expression for $E(\Omega|S)$ and $E(\Sigma|S)$ when the prior distribution on $\sigma = 2\Sigma$ is the objective prior (4.4). From Proposition 3.2 and (3.9), we immediately obtain the following.

COROLLARY 4.1. *Let $Z_i, i = 1, \ldots, n$ be a sample from the $N(0, \Sigma)$ distribution with $\Sigma \in Q_G$. Let $U = \sum_{i=1}^n Z_i Z_i^t$ and let the prior distribution on $2\Sigma$ be as in (4.4). Then the posterior mean of $\Omega = \widehat{\Sigma}^{-1}$ is*

$$
\begin{aligned}
E(\Omega|S) = \sum_{j=1}^k (S_{C_j}^{-1})^0 &- \left(1 - \frac{c_1 + c_2 - 2s_2}{n}\right)(S_{S_2}^{-1})^0 \\
&- \sum_{j=3}^k \left(1 - \frac{c_j - s_j}{n}\right)(S_{S_j}^{-1})^0.
\end{aligned}
$$

(4.7)

It is interesting to note here that when $n$ tends to $+\infty$, the expression of the posterior mean in (4.7) becomes very close to the expression of the $mle_g$ of $\Omega$ as given in (2.9). This will also be illustrated by our numerical results further in this paper.

From Theorem 3.1, we immediately derive the posterior mean of $\Sigma$ as follows.

COROLLARY 4.2. *Let $U$ and $\Sigma$ as above; then $E(2\Sigma|S)$ is given by (3.10)–(3.14) for $\theta = \pi(nS)$ and $(\alpha, \beta)$ as in (4.6) with the additional condition that $(\alpha, \beta)$ in (4.6) satisfy the inequalities*

$$
\alpha_1 + \frac{c_1 + 1}{2} + \gamma_2 < 0, \qquad \alpha_1 + \frac{c_1 - s_2 + 1}{2} < 0,
$$

$$
\alpha_j + \frac{c_j - s_j + 1}{2} < 0, \qquad j = 2, \ldots, k.
$$

We note here that the additional conditions imposed on the posterior hyperparameters are there to insure that the moments given in (3.10)–(3.14) exist.

**5. Decision-theoretic results.** In this section, we will derive the Bayes estimators for $\Sigma$ and $\Omega$, under two loss functions similar to the classical $L_1$ and $L_2$ loss functions, but adapted to $Q_G$ and $P_G$. These Bayes estimators are of course computed with respect to a given prior distribution. The prior distributions we will consider are the $HIW$ and the $IW_{P_G}$ as recalled in Section 3, and the reference prior as developed in Section 4.

5.1. *Bayes estimators for $\Sigma$ and $\Omega$.* We now proceed to place the covariance estimation problem in graphical models in a decision-theoretic framework. Let us first recall what is traditionally done in the saturated case, that is, when $G$ is complete. Given a sample of size $n$ from a $N_r(0, \Sigma)$ distribution, letting $\tilde{\Sigma}$ be any estimator of $\Sigma$ based on that sample, we consider the following two loss functions:

$$
\begin{aligned}
L_1(\tilde{\Sigma}, \Sigma) &= \langle \tilde{\Sigma}, \Sigma^{-1} \rangle - \log |\tilde{\Sigma}\Sigma^{-1}| - r, \\
L_2(\tilde{\Sigma}, \Sigma) &= \langle \tilde{\Sigma} - \Sigma, \tilde{\Sigma} - \Sigma \rangle,
\end{aligned}
\tag{5.1}
$$

called Stein's (or entropy, or likelihood [19, 34]) and squared-error (or Frobenius [25, 28]) losses, respectively. Other losses have also been considered in the literature (see [33] for details). Many authors such as Haff [15], Krishnamoorthy [23] and Krishnamoorthy and Gupta [22] have also considered the estimation of the precision matrix $\Omega = \Sigma^{-1}$, instead of $\Sigma$. The reader is referred to [37] for a more complete list. The natural analogues for $\Omega$ of (5.1) are

$$
\begin{aligned}
L_1(\tilde{\Omega}, \Omega) &= \langle \tilde{\Omega}, \Omega^{-1} \rangle - \log |\tilde{\Omega}\Omega^{-1}| - r, \\
L_2(\tilde{\Omega}, \Omega) &= \langle \tilde{\Omega} - \Omega, \tilde{\Omega} - \Omega \rangle = \mathrm{tr}(\tilde{\Omega} - \Omega)^2.
\end{aligned}
\tag{5.2}
$$

A question that naturally arises in various contexts in multivariate analysis and related topics is whether to estimate $\Sigma$ or its inverse $\Omega = \Sigma^{-1}$. We choose to focus on the estimation of both $\Sigma$ and $\Omega$ in this paper for a variety of reasons. The parameter $\Sigma$ has a natural and well-understood interpretation in multivariate analysis and its direct estimation has numerous applications. The precision matrix, on the other hand, has a natural and central place in Gaussian graphical models as it is the canonical parameter of the natural exponential family $\mathcal{N}_G$ and it sits in the parameter set $P_G$ as defined in (2.2), a parameter set of dimension much smaller than that of $M_r^+$.

To our knowledge, in the case where $G$ is decomposable and not complete, a decision-theoretic estimation of the scale parameters $\Sigma$ or $\Omega$ has not been previously considered. We will do so now. We first observe that the traditional loss functions used for saturated models need to be reconsidered for graphical models as we now have fixed zeros in the inverse covariance matrix and therefore fewer parameters. This is clear in the expression of $L_2(\tilde{\Omega}, \Omega)$ in (5.2) when $\Omega$ is in $P_G$. Indeed, we have

$$
L_2(\tilde{\Omega}, \Omega) = \sum_{(i,j) \in E} (\tilde{\Omega}_{ij} - \Omega_{ij})^2
$$

dependent not on $r(r+1)/2$ parameters but on the nonzero parameters only. Similarly, because of the structural zeros of $\Omega \in P_G$, $L_1(\tilde{\Omega}, \Omega)$ in (5.2) depends only on the nonzero elements of $\Omega$. We also note that, according to (2.7) and (2.8), for an arbitrary decomposable graph $G$, when $\widetilde{\Sigma}$ and $\Sigma$ both belong to $Q_G$, (5.1) can be written as $L_1(\widehat{\widetilde{\Sigma}}, \widehat{\Sigma}) = L_1(\widetilde{\Sigma}, \Sigma)$, where

$$
\begin{aligned}
L_1(\widetilde{\Sigma}, \Sigma) = \sum_{C \in \mathcal{C}} \langle \widetilde{\Sigma}_C, \Sigma_C^{-1} \rangle - \sum_{S \in \mathcal{S}} \langle \widetilde{\Sigma}_S, \Sigma_S^{-1} \rangle \\
- \log \frac{\prod_{C \in \mathcal{C}} |\widetilde{\Sigma}_C| \prod_{S \in \mathcal{S}} |\widetilde{\Sigma}_S|}{\prod_{C \in \mathcal{C}} |\Sigma_C| \prod_{S \in \mathcal{S}} |\widetilde{\Sigma}_S|} - r,
\end{aligned}
$$

(5.3)

which involves solely the elements of $\Sigma \in Q_G$ and $\widetilde{\Sigma} \in Q_G$ and not the nonfree elements of their completions $\widehat{\Sigma}$ and $\widehat{\widetilde{\Sigma}}$. Accordingly, we shall modify the traditional $L_2$ loss function for $\Sigma \in Q_G$ as follows so that, like $L_1(\widetilde{\Sigma}, \Sigma)$, it depends only on the free parameters of $\widehat{\widetilde{\Sigma}}$. For $\Sigma$ and $\widetilde{\Sigma}$ in $Q_G$, we define

$$
(5.4) \qquad L_2(\widetilde{\Sigma}, \Sigma) = \langle \widetilde{\Sigma} - \Sigma, \widetilde{\Sigma} - \Sigma \rangle = \sum_{(i,j) \in E} (\widetilde{\Sigma}_{ij} - \Sigma_{ij})^2.
$$

We have therefore modified, when necessary, the traditional loss functions given by (5.1)–(5.2) to take into account the graphical nature of the covariance matrix. We now derive the corresponding Bayes estimators for these newly defined loss functions and given priors.

PROPOSITION 5.1. *Let $Z_i, i = 1, \ldots, n$ and $U = nS$ be as in Corollary 3.1. Then, for a given prior $\pi(\Sigma)$ on $\Sigma \in Q_G$, the Bayes estimators of $\Sigma$ under (5.3) and (5.4) are equal to, respectively,*

$$
(5.5) \qquad \widetilde{\Sigma}_{L1}^{\pi(\Sigma|U)} = \kappa([E^{\pi(\Sigma|U)}(\widehat{\Sigma}^{-1})]^{-1}) \quad and \quad \widetilde{\Sigma}_{L2}^{\pi(\Sigma|U)} = E^{\pi(\Sigma|U)}(\Sigma),
$$

*where $\pi(\Sigma|U)$ denotes the posterior distribution of $\Sigma$ given $U$.*

*For a given prior $\pi(\Omega)$ for $\Omega \in P_G$, the Bayes estimators of $\Omega$ under the loss functions in (5.2) are equal to, respectively,*

$$
(5.6) \qquad \widetilde{\Omega}_{L1}^{\pi(\Omega|U)} = [E^{\pi(\Omega|U)}\widehat{(\kappa(\Omega^{-1}))}]^{-1} \quad and \quad \widetilde{\Omega}_{L2}^{\pi(\Omega|U)} = E^{\pi(\Omega|U)}(\Omega).
$$

PROOF. We first derive the expression of $\widetilde{\Sigma}_{L1}^{\pi(\Sigma|U)}$ in (5.5). The Bayes estimator is the estimator $\widetilde{\Sigma}$ that minimizes the posterior expected loss. So for $L_1$ loss, we have

$$
\begin{aligned}
E^{\pi(\Sigma|U)}[L_1(\widetilde{\Sigma}, \Sigma)] &= \int [\langle \widetilde{\Sigma}, \widehat{\Sigma}^{-1} \rangle - \log |\widetilde{\Sigma}\widehat{\Sigma}^{-1}| - r] \pi(\Sigma|U) \, d\Sigma \\
&= \langle \widetilde{\Sigma}, E^{\pi(\Sigma|U)}[\widehat{\Sigma}^{-1}] \rangle - \log |\widetilde{\Sigma}| - E^{\pi(\Sigma|U)}[\log |\widehat{\Sigma}^{-1}|] - r \\
&= \langle \widetilde{\Sigma}, R(U) \rangle - \log |\widetilde{\Sigma}| - c
\end{aligned}
$$

for some constant $c$ and $R(U) = E^{\pi(\Sigma|U)}[\widehat{\Sigma}^{-1}]$. Minimizing the posterior expected loss is equivalent to maximizing the function $\log |\widetilde{\Sigma}| - \langle \widetilde{\Sigma}, R(U) \rangle$ with respect to $\widetilde{\Sigma}$. This function is concave and is maximized at $\widetilde{\Sigma} = R(U)^{-1} = [E^{\pi(\Sigma|U)}[\widehat{\Sigma}^{-1}]]^{-1}$ hence yielding the Bayes estimator of $\Sigma$ for $L_1$ loss.

Let us now derive the expression of $\widetilde{\Sigma}_{L2}^{\pi(\Sigma|U)}$ in (5.5). Once more the $L_2$ Bayes estimator is found by minimizing the posterior expected loss. Now,

$$R_{L_2}[\widetilde{\Sigma}, \Sigma] = \int \sum_{(i,j) \in E} (\widetilde{\Sigma}_{ij} - \Sigma_{ij})^2 \pi(\Sigma|U) \, d\Sigma$$

$$= \sum_{(i,j) \in E} \int (\widetilde{\Sigma}_{ij} - \Sigma_{ij})^2 \pi(\Sigma|U) \, d\Sigma.$$

Since we are minimizing a sum of terms, it is sufficient to minimize each one of the terms. It is well known (see [15]) that each one of these is minimized for $\tilde{\Sigma}_{ij} = E^{\pi(\Sigma|U)}(\Sigma_{ij})$, which gives us the desired expression for $\tilde{\Sigma}_{L_2}^{\pi(\Sigma|U)}$. The proofs for the Bayes estimators of $\Omega$ follow along similar lines. □

We derived the expression of $\widetilde{\Sigma}_{L2}^{\pi(\Sigma|U)}$ above, even though it is straightforward, in order to emphasize the fact that, unlike in the classical case when one estimates a complete covariance matrix, the posterior means are Bayes optimal only if the $L_2$ loss function is modified to reflect the graph that underlies $\Sigma$. In other words the posterior mean will not be the Bayes estimator if nonfree elements of the matrix $\Sigma$ contribute to the loss, a point that can be easily overlooked when considering decision-theoretic estimation for graphical models using the traditional loss functions.

It is important to note here, since it will simplify many of our computations in Sections 6 and 7, that the $L_1$ estimator for $\Sigma$ is the $\varphi$ transformation of the $L_2$ estimator for $\Omega$. A similar result holds for the $L_2$ estimator for $\Sigma$, that is,

$$\tilde{\Sigma}_{L_1} = \kappa([\tilde{\Omega}_{L_2}]^{-1}), \qquad \tilde{\Sigma}_{L_2} = \kappa([\tilde{\Omega}_{L_1}]^{-1}).$$

The risk functions corresponding to the losses above are

$$R_{L_i}(\tilde{\Sigma}_{L_i}) = E[L_i(\tilde{\Sigma}_{L_i}, \Sigma)], \qquad R_{L_i}(\tilde{\Omega}_{L_i}) = E[L_i(\tilde{\Omega}_{L_i}, \Omega)], \qquad i = 1, 2.$$

In the subsequent sections, these risk functions will be used to assess the quality of the eight estimators that we consider. For each of $\Sigma$ and $\Omega$, the eight estimators considered will be the sample covariance matrix $S$ (if ignoring the graphical model structure) and its inverse, the $mle_g$ for $\Sigma$ and $\Omega$, $\widetilde{\Sigma}_g$ and $\widetilde{\Omega}_g$, and

$$\widetilde{\Sigma}_{L_i}^{\pi}, \qquad i = 1, 2 \quad \text{and} \quad \widetilde{\Omega}_{L_i}^{\pi}, \qquad i = 1, 2,$$

where the prior $\pi$ will be a $HIW$ or more generally an $IW_{P_G}$, or the reference prior.

5.2. *Risk properties of* $\widetilde{\Sigma}_g$. We now proceed to state a decision-theoretic property of the maximum likelihood estimator $\widetilde{\Sigma}_g$ of $\Sigma \in Q_G$. From (2.9), it follows immediately that

$$\widetilde{\Sigma}_g = \left[ \sum_{j=1}^k [S_{C_j}^{-1}]^0 - \sum_{j=2}^k [S_{S_j}^{-1}]^0 \right]^{-1}.$$

LEMMA 5.1. *The maximum likelihood estimator* $\widetilde{\Sigma}_g$ *is the best* $L_1$ *estimator in the class of estimators of the form* $a\widetilde{\Sigma}, a \in \mathbb{R}$ *where*

$$(5.7) \qquad \widetilde{\Sigma} = \left[ \sum_{j=1}^k [\kappa(U)_{C_j}^{-1}]^0 - \sum_{j=2}^k [\kappa(U)_{S_j}^{-1}]^0 \right]^{-1}$$

*and* $U = \sum_{i=1}^n Z_i Z_i^t$.

PROOF. Recall that for the estimator $a\widetilde{\Sigma}$ under $L_1$ loss we have

$$R_1(a\widetilde{\Sigma}, \Sigma) = E[\langle \widehat{\Sigma}^{-1}, a\widetilde{\Sigma} \rangle - \log |\widehat{\Sigma}^{-1} a\widetilde{\Sigma}| - r]$$
$$= a\langle \widehat{\Sigma}^{-1}, E[\widetilde{\Sigma}] \rangle - r\log(a) - E[\log |\widehat{\Sigma}^{-1}\widetilde{\Sigma}|] - r.$$

Now since $E[\widetilde{\Sigma}] = n\Sigma$ (see [24], page 133) we have that $a\langle \widehat{\Sigma}^{-1}, E[\widetilde{\Sigma}] \rangle = nra$. Moreover, by (2.8),

$$E[\log |\widehat{\Sigma}^{-1}\widetilde{\Sigma}|] = E\left[ \log \frac{\det \widetilde{\Sigma}}{\det \Sigma} \right] = E\left[ \log \frac{\prod_{C \in \mathcal{C}} |\widetilde{\Sigma}_C| \prod_{S \in \mathcal{S}} |\Sigma_S|}{\prod_{C \in \mathcal{C}} |\Sigma_C| \prod_{S \in \mathcal{S}} |\widetilde{\Sigma}_S|} \right]$$
$$= E\left[ \log \frac{\prod_{C \in \mathcal{C}} |\Sigma_C^{-(1/2)} \widetilde{\Sigma}_C \Sigma_C^{-(1/2)}|}{\prod_{S \in \mathcal{S}} |\Sigma_S^{-(1/2)} \widetilde{\Sigma}_S \Sigma_S^{-(1/2)}|} \right].$$

Since $\Sigma_C^{-(1/2)} \widetilde{\Sigma}_C \Sigma_C^{-(1/2)} \sim W_C(n, I_C)$ and $\Sigma_S^{-(1/2)} \widetilde{\Sigma}_S \Sigma_S^{-(1/2)} \sim W_S(n, I_S)$ (see [31]),

$$|\Sigma_C^{-(1/2)} \widetilde{\Sigma}_C \Sigma_C^{-(1/2)}| = \prod_{j=1}^{c_i} \chi_{n-j+1}^2, \qquad |\Sigma_S^{-(1/2)} \widetilde{\Sigma}_S \Sigma_S^{-(1/2)}| = \prod_{j=1}^{s_i} \chi_{n-j+1}^2.$$

Letting $E[\log |\widehat{\Sigma}^{-1}\widetilde{\Sigma}|] = m$ which is a constant independent of $\Sigma$ or $a$, we therefore now have $R_1(a\widetilde{\Sigma}, \Sigma) = anr - r\log a - m - r$. Differentiating with respect to $a$ and setting the derivative to zero gives $a = 1/n$, which proves the lemma. $\square$

**6. Risk comparisons and numerical properties.** In this section, through two examples, we investigate the performance of our Bayes estimators derived from the different priors presented in Sections 3 and 4. We base our comparisons on frequentist risk calculations obtained from simulations under losses $L_1$ and $L_2$ for both $\Sigma$ and $\Omega$, on predictive properties and on eigenvalues properties.

6.1. *Example 1: "Two Cliques."* In this example, we illustrate the power and flexibility of the $IW_{P_G}$ family and its multiple shape parameters. We build a general example based on the call center data analyzed in [18] and described in Section 7. First, we define a graph $G$ with 100 vertices ($r = 100$) where $C_1 = \{1, \ldots, 70\}$, $C_2 = \{61, \ldots, 100\}$ and $S_2 = \{61, \ldots, 70\}$. The true covariance matrix $\Sigma$ is constructed by (i) taking the sample covariance of the first 100 variables in the call center data, (ii) removing the $ij$ entries corresponding to $(i, j)$ which are not edges of $G$ and (iii) performing the completion operation described in (2.1). This procedure guarantees that $\Sigma$ preserves the conditional independence relationships specified in $G$. This example involves cliques of different dimensions ($c_1 = 70$ and $c_2 = 40$) and our goal is to show that it is possible to obtain improved estimators by using an $IW_{P_G}$ prior with multiple shape parameters. The idea here is to apply different levels of shrinkage for each clique, that is, more shrinkage for larger cliques, following the intuition that more shrinkage is necessary in higher-dimensional problems.

Our simulations compare estimators based on the reference prior, the traditional hyper inverse Wishart (one shape parameter, $\delta = 3$) and five versions of the $IW_{P_G}$. In the latter, we first choose the shape parameters in proportion to clique size. More specifically, we choose $\alpha_i$ of the form

$$\alpha_i = -\frac{\delta_i + c_i - 1}{2}, \qquad i = 1, 2,$$

with $\delta_i$ equal to $\frac{1}{2}c_i$, $\frac{1}{4}c_i$ and $\frac{1}{10}c_i$. In addition, we also compare risk obtained from empirical Bayes estimates of shape parameters. These were computed based on the following specifications of $\alpha$ and $\beta$ as a function of clique size:

(i) $\alpha_i = \frac{-(\delta c_i + c_i - 1)}{2}$ and $\beta_i = \frac{-(\delta s_i + s_i - 1)}{2}$,

(ii) $\alpha_i = ac_i + b$ and $\beta_i = as_i + b$.

Notice that in order to obtain the empirical Bayes estimates of $\alpha$ and $\beta$ we maximize the marginal likelihood as a function of $\delta$ in case (i) whereas in (ii) the maximization is done as a function of both $a$ and $b$. For all choices of shape parameters, we define the scale matrix in two ways: the identity matrix $I$ and a matrix $D$ for which the prior expected value of $\Sigma \in Q_G$ is $I$. Conditional on a shape parameter, $D$ can be easily derived as we saw in Lemma 3.1. This example focuses on $L_1$ loss and the impact of flexible priors in estimating the eigenstructure of $\Sigma$ and $\Omega$. The results in Table 1 show that appropriate choices of different shape parameters have a significant effect in reducing the risk of estimators. Looking at $L_1(\Omega)$, and $n = 75$, for the estimator under $IW_{P_G}(\frac{1}{4}c_i, D)$ there is approximately a 78% (76% for $n = 100$, 20% for $n = 500$ and 9% for $n = 1000$) reduction in risk when compared with the more traditional $HIW(3, I)$. When comparing to the constrained maximum likelihood estimator ($mle_g$) the reduction is even more impressive: 97% for $n = 75$, 87% for $n = 100$, 30% for $n = 500$ and 16% for $n = 1000$. Additionally the reference prior performs well and always beats the $mle_g$.

TABLE 1
*Two cliques: risk for estimators*

| | $n = 75$ | | $n = 100$ | | $n = 500$ | | $n = 1000$ | |
|---|---|---|---|---|---|---|---|---|
| | $L_1(\Omega)$ | $L_1(\Sigma)$ | $L_1(\Omega)$ | $L_1(\Sigma)$ | $L_1(\Omega)$ | $L_1(\Sigma)$ | $L_1(\Omega)$ | $L_1(\Sigma)$ |
| *Reference* | 212.7 | 66.61 | 60.71 | 40.93 | 7.02 | 6.66 | 3.33 | 3.28 |
| $HIW(3, I)$ | 98.76 | 59.28 | 80.72 | 43.41 | 7.76 | 7.18 | 3.54 | 3.43 |
| *Empirical Bayes* (i) | 127.47 | 58.43 | 78.05 | 42.87 | 7.70 | 7.13 | 3.52 | 3.41 |
| *Empirical Bayes* (i)-D | 27.44 | 25.84 | 23.74 | 21.95 | 6.31 | 6.02 | 3.20 | 3.12 |
| *Empirical Bayes* (ii) | 121.55 | 57.12 | 74.45 | 41.81 | 7.60 | 7.05 | 3.51 | 3.39 |
| *Empirical Bayes* (ii)-D | 24.78 | 23.04 | 21.48 | 19.95 | 6.12 | 5.86 | 3.16 | 3.08 |
| $IW_{P_G}(1/2c_i, D)$ | 29.99 | 25.18 | 24.53 | 24.49 | 6.37 | 6.21 | 3.27 | 3.22 |
| $IW_{P_G}(1/2c_i, I)$ | 207.4 | 67.88 | 116.7 | 49.78 | 8.69 | 7.80 | 3.76 | 3.61 |
| $IW_{P_G}(1/4c_i, D)$ | 18.57 | 15.87 | 18.57 | 15.87 | 5.67 | 5.43 | 3.03 | 2.96 |
| $IW_{P_G}(1/4c_i, I)$ | 165.5 | 63.10 | 96.14 | 46.20 | 8.14 | 7.43 | 3.67 | 3.50 |
| $IW_{P_G}(1/10c_i, D)$ | 35.71 | 31.99 | 31.59 | 27.02 | 6.77 | 6.41 | 3.32 | 3.23 |
| $IW_{P_G}(1/10c_i, I)$ | 141.7 | 60.23 | 89.67 | 45.03 | 7.98 | 7.32 | 3.59 | 3.47 |
| *MLEg* | 813.9 | 70.72 | 154.6 | 43.51 | 8.13 | 6.79 | 3.62 | 3.32 |
| *MLE* | – | – | $7.3 \times 10^8$ | 102.5 | 14.45 | 10.85 | 6.00 | 5.22 |

We emphasize the connection between $L_1$ and the eigenstructure of $\Omega$ and $\Sigma$ by the scree plots in Figure 1. It is our belief that the superior performance of the Bayesian estimators under the $IW_{P_G}$ and the reference prior relative to the $mle_g$ is a direct consequence of the better estimation of the eigenvalues of both the precision and covariance matrices.

This experiment also indicates that choosing the amount of shrinkage as a function of shape parameters can be a delicate task. Here, we find that $\delta_i = \frac{1}{4}c_i$ (i.e., $\delta_1 = 17.5$ and $\delta_2 = 10$) performs best and seems to be a good compromise between $\delta_i = \frac{1}{10}c_i$ and $\delta_i = \frac{1}{2}c_i$. The definition of appropriate hyperparameters and the choice of shrinkage level is context specific and depends on the amount of prior information available. An alternative to the subjective specification of the hyperparameters is the empirical Bayes approach presented. The results in Table 1 show that this alternative performs reasonably well and uniformly outperforms the reference prior (the other objective alternative). As we show in this example, the $IW_{P_G}$ offers a very general framework for the incorporation of prior knowledge with the ability to significantly improve the performance of posterior estimators and estimation of covariance matrices in general.

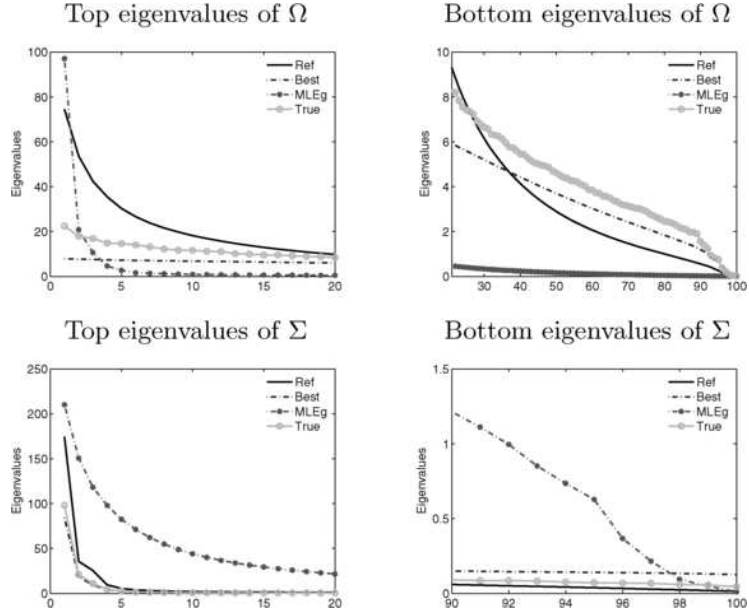

FIG. 1. *Scree plots of eigenvalues for $\Omega$ (top row) and $\Sigma$ (bottom row) in the simulation with $n = 100$. For each estimator, the lines represent the average of the eigenvalues after 1000 simulations. "Best" refers to the estimator with lowest risk $(IW_{P_G}(1/4c_i, D))$.*

6.2. *Example 3: choosing the graph.* Our second example demonstrates the potential of graphical models as a model-based tool for regularization and estimation of large covariance structures. So far, we have presented examples that compare the performance of different estimators (based on different priors) assuming knowledge of $G$. In real problems, $G$ is unknown and often has to be inferred before parameter estimation. From a Bayesian perspective, model selection involves the exploration of the posterior distribution of graphs given by

$$(6.1) \qquad p(G|X) \propto p(X|G)p(G),$$

where $p(X|G)$ is the marginal likelihood of $G$ and $p(G)$ represents its prior. In Gaussian graphical models the marginal likelihood for any $G$ is given by the following integral:

$$(6.2) \qquad p(X|G) = \int_{Q_G} f(X|\Sigma, G)\pi(\Sigma|G)\,d\Sigma,$$

where $f(X|\Sigma, G)$ is the density of $X$ given $\Sigma$ and $G$ and $\pi(\Sigma|G)$ is the prior distribution for $\Sigma$ given $G$. Using $IW_{P_G}$ conjugate priors for $\Sigma$ makes the computation of the above integral straightforward since the expression of the marginal likelihood is obtained explicitly through the normalizing constant of the $IW_{P_G}$ as given in (3.4). If we assume a uniform prior over the graphs, computing the posterior distribution of graphs is equivalent to computing the marginal likelihoods.

To illustrate how graphical models can be used as a regularization tool, we build an example where the underlying graph is unknown and will be selected based on marginal likelihoods of a particular restricted set of graphs. We focus on a subclass $\mathcal{G}_{b60}$ of graphs where the precision matrix is "banded" (see [2]). The restricted subclass $\mathcal{G}_{b60}$ consists of the decomposable graphs $G_k, k = 1, \ldots, 60$ with cliques $C_j = \{j, j+1, \ldots, j+k\}, j = 1, \ldots, r-k$. The graphical Gaussian model Markov with respect to $G_k$ can be viewed as an $AR(k)$ model and the corresponding precision matrix is a banded matrix with a band of width $k + 1$, indicating that all elements beyond the $k$th supra-diagonals are zero. For added simplicity, we use the $HIW(3, I)$, a special case of the $IW_{P_G}$, as a prior for $\Sigma$.

As in the previous example we build the true covariance matrix from the sample covariance of the call center data using a graph corresponding to $k = 20$ followed by the completion operation in (2.1). We proceed by sampling $n$ observations from a $N_r(0, \Sigma)$, 1000 times. At each iteration, we

TABLE 2
*Choosing the graph example: risk for estimators*

|  | $R_1(\Omega)$ | $R_1(\Sigma)$ | $R_2(\Omega)$ | $R_2(\Sigma)$ |
|---|---|---|---|---|
|  |  | $n = 100$ | $\hat{k} = 4.36$ |  |
| *Reference* | 15.36 (23.67) | 17.14 (19.85) | 1902.6 (13375.0) | 22.19 (241.2) |
| $HIW(3, I)$ | 14.76 (22.77) | 16.55 (19.11) | 1736.3 (10350.0) | 16.23 (56.88) |
| *MLEg* | 15.89 (32.64) | 18.42 (21.08) | 1876.0 (9897.80) | 16.54 (57.74) |
| *MLE* | $9.9 \times 10^6$ | 102.53 | $1.1 \times 10^{18}$ | 133.08 |
|  |  | $n = 500$ | $\hat{k} = 7.38$ |  |
| *Reference* | 8.084 (11.94) | 9.078 (11.79) | 1105.6 (1500.6) | 5.648 (16.01) |
| $HIW(3, I)$ | 8.053 (11.96) | 8.961 (12.81) | 1101.7 (1234.5) | 5.070 (11.09) |
| *MLEg* | 8.129 (12.20) | 9.297 (15.85) | 1116.8 (1571.2) | 5.088 (11.12) |
| *MLE* | 14.45 | 10.85 | 3147.4 | 26.43 |
|  |  | $n = 1000$ | $\hat{k} = 18.80$ |  |
| *Reference* | 2.256 (1.930) | 2.203 (1.893) | 345.5 (310.3) | 6.624 (7.009) |
| $HIW(3, I)$ | 2.259 (1.936) | 2.197 (1.899) | 350.4 (317.9) | 5.480 (5.736) |
| *MLEg* | 2.311 (1.992) | 2.232 (1.910) | 331.8 (291.3) | 5.492 (5.747) |
| *MLE* | 6.003 | 5.226 | 1006.8 | 13.20 |

The values in parentheses refer to the risk of estimators generated by the "oracle" $(k = 20)$.

compute the marginal likelihood for all graphs with $k = 1, \ldots, 60$, choosing the top model to proceed with the estimation of $\Sigma$. For each estimator, under the different priors, losses are computed and risk compared. We also investigate the performance of the procedure by comparing estimators generated through the "oracle" that knows the correct value of $k$. The corresponding risk values are given in brackets.

We repeat this exercise for $n = 100, 500$ and $1000$ with results presented in Table 2. In each case, the average of the best $k$'s is denoted $\hat{k}$ and is also given in the table.

It is clear from the example that more parsimonious models are selected from small sample sizes. Indeed, for $n = 100$ we choose the average band size $\hat{k} = 4.36$; for $n = 500$, we choose $\hat{k} = 7.38$ and for $n = 1000$, $\hat{k} = 18.80$. This highlights the *Ockham's razor* effect of marginal likelihoods [20], in selecting a graph. Moreover, for $n = 100$ and $n = 500$ the losses generated by the oracle are always larger than those of our estimators, with the oracle only being relatively competitive for $n = 1000$.

**7. Call center data.** In this section, we apply our methodology to the call center data analyzed in [2, 18]. With this example, we will illustrate the predictive properties of our estimators and the flexibility yielded by graphical models. Indeed, we will show that our estimators when using banded matrices and the $IW_{P_G}$ have a smaller predictive error than the $mle_g$. More strikingly, we will show that when using bands varying in width along the diagonal together with the $IW_{P_G}$, our estimators yield significantly improved predictive power over the best uniformly banded model.

The dataset in this example constitutes records from a call center of a major financial institution in 2002 where the number of incoming phone calls during 10-minute intervals from 7:00 am till midnight were recorded. Weekends, public holidays and days with equipment malfunction are excluded, resulting in data for 239 days. The number of calls in each of these intervals is denoted as $N_{ij}$, $i = 1, 2, \ldots, 239$ and $j = 1, 2, \ldots, 102$. A standard transformation $x_{ij} = (N_{ij} + \frac{1}{4})^{1/2}$ is applied to the raw data to make it closer to Normal.

7.1. *Analysis with banded precision matrices.* We consider the class of models $\mathcal{G}_{b60}$ as described in Section 6.2 and as we saw there, choosing a model Markov with respect to $G_k \in \mathcal{G}_{b60}$ is equivalent to banding the inverse covariance matrix and our results can therefore be readily compared to those of Bickel and Levina [2]. It is important to note, however, that our approach differs from that of [2] in two ways. First, in [2] the banded estimators for the precision matrix are used as estimators of the saturated $\Sigma^{-1}$. We, on the other hand, fit graphical models to the call center data, explicitly assuming that the true $\Sigma$ is such that $\Sigma^{-1}$ is in $P_G$, that is, has fixed zeros. Second, we use the eight estimators for $\Sigma$ as described in Section 5.1. This includes the traditional frequentist estimator for the graphical model, $mle_g$, and the Bayesian estimators that we have developed above.

We employ both the traditional cross-validation and Bayesian model selection procedures to determine the "best" model among the class of graphical models with $k$-banded precision matrices. The cross-validation procedure is done through the $K$-fold cross-validation method with $K = 10$. The dataset with 239 data points is divided into 10 parts (the first nine parts have 40 observations and the last has 39 observations). We predict the second half of the day, given data for the first half, on the test data after computing estimators on the training dataset. In particular, we partition the 102-dimensional random vectors $X_1, X_2, \ldots, X_{239}$ into two equal parts, each representing the first and second half of the day as follows:

$$x = \begin{pmatrix} x^{(1)} \\ x^{(2)} \end{pmatrix}, \qquad \mu = \begin{pmatrix} \mu_1 \\ \mu_2 \end{pmatrix}, \qquad \Sigma = \begin{pmatrix} \Sigma_{11} & \Sigma_{12} \\ \Sigma_{21} & \Sigma_{22} \end{pmatrix},$$

where $x_i^{(1)} = (x_{i1}, x_{i2}, \ldots, x_{i51})^t$ and $x_i^{(2)} = (x_{i52}, x_{i53}, \ldots, x_{i102})^t$. The mean vectors are partitioned in a similar manner. The best linear predictor for $x_i^{(2)}$ from $x_i^{(1)}$ is given by

$$x_i^{(2)} = \mu_2 + \Sigma_{21} \Sigma_{11}^{-1} (x_i^{(1)} - \mu_1).$$

We use the prediction equation above with the following eight estimators for $\Sigma$: the $mle$ and $mle_g$, the estimators based on $L_1$ loss, $\widetilde{\Sigma}_{L_1}^{IW_{P_G}}, \widetilde{\Sigma}_{L_1}^{HIW}, \widetilde{\Sigma}_{L_1}^{Ref}$, and the estimators based on $L_2$ loss, $\widetilde{\Sigma}_{L_2}^{IW_{P_G}}, \widetilde{\Sigma}_{L_2}^{HIW}, \widetilde{\Sigma}_{L_2}^{Ref}$. For the Bayes estimators we use the traditional choice of the identity matrix as the scale hyperparameter and the shape parameters for the $IW_{P_G}$ are set as $\alpha_i = -5, \forall i, \beta_2 = -4 + \frac{k}{2}$. The prediction error on the test data is measured by the average absolute forecast error. The Bayesian model selection procedure entails choosing the model with the maximum marginal likelihood according to the principles given in Section 6.2. Here the prior distribution used for $\Sigma \in Q_G$ is the $IW_{P_G}$ with the same hyperparameters as above.

For all eight estimators, the cross-validation procedure identifies $k = 4$ as the model with the lowest prediction error for the second half of the day given the first. The Bayesian model selection procedure identifies $k = 5$ as the model with highest marginal likelihood. We note that both model selection procedures yield very similar, parsimonious models for the call center data.

We proceed to compare the forecast performance (or prediction error) of our estimators. As done in [2], we forecast the second half of the day based on the first half using the first 205 data points as the training set and the last 34 as the test set. The prediction error lines are given in Figure 2.

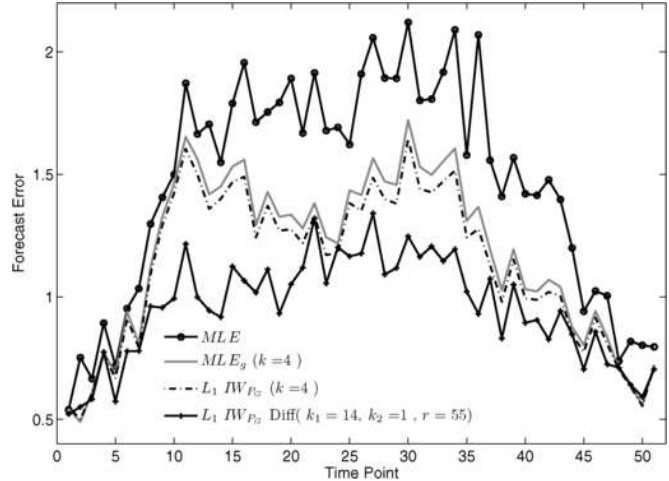

Fig. 2. *Forecast error for selected banded and "differentially banded" models.*

The $\widetilde{\Sigma}_{L1}^{IW_{P_G}}$ and $\widetilde{\Sigma}_{L1}^{Ref}$ prediction error lines are so close that it is difficult to distinguish between the two in a plot and thus we only show the one for $\widetilde{\Sigma}_{L1}^{IW_{P_G}}$. Prediction errors using the Bayes estimators are all lower than the prediction error given by the $mle_g$ and by the standard sample covariance matrix, the $mle$. For the sake of clarity, Figure 2 gives the forecast error for the best model ($AR(4)$) chosen via cross-validation for the $mle$, the $mle_g$, $\widetilde{\Sigma}_{L1}^{IW_{P_G}}$ and another estimator, "$L_1 IW_{P_G}$diff" we will describe in Section 7.2.

Besides the overall poor performance of the standard estimator $S$, it is also well-understood that $S$ overestimates the largest eigenvalues and underestimates the lowest eigenvalue. An examination of the scree plots implied by the eigenvalues of our Bayes estimators (not shown here for brevity) reveals that our Bayes estimators compared to $S$ have lower estimates of the largest eigenvalues and higher estimates for the smallest eigenvalues—hence the Bayes estimators seem to shrink the eigenvalues toward the center of the eigenspectrum, a property often sought by estimators proposed in previous work [9, 25, 35].

7.2. *Analysis with differential banding.* In the above analysis we restricted ourselves to the class of $k$-banded inverse covariance models or $AR(k)$ models. This approach highlighted, among other important properties, the fact that banding the inverse covariance matrix, as carried out in Bickel and Levina [2], essentially entails fitting graphical models.

We noted that the cross-validation error from predicting the second half of the day given the first half suggested that the $AR(4)$ model (i.e., $k = 4$) gave the lowest prediction error. The cross-validation error from predicting the first half of the day given the second half suggested that the $AR(16)$ model (i.e., $k = 16$) gives the lowest prediction error. These two different values of $k$ suggests that one single $k$ may not be sufficient to explain the features of this dataset. An examination of the sample correlation matrix suggests different correlation structures during the first and second half of the day. In particular the correlation between variables which are farther apart is stronger in the first half of the day than in the second half of the day, suggesting that the order of the lag, which represents the level of the strength between neighboring variables, could be different in different parts of the day. It seems that "differential banding" of the inverse covariance matrix is necessary to capture this effect and fitting a straight $k$-banded model, as was done previously, may not necessarily allow the flexibility that we want.

A natural approach to obtaining this flexibility is to frame the problem once more in the context of graphical models. Let us now consider graphical models with two different clique sizes for the two parts of the day as an extension of the single clique size suggested by $k$-banded models. We note that $k$-banded models have cliques of size $k + 1$ and separators of size $k$. Let us consider what we term as $(k_1, k_2, r)$ "differentially banded" models, where $k_1 + 1$ represents the size of the cliques in the first part of the day, $k_2 + 1$ the size of the cliques in the second part of the day, and $r$ the point at which the change takes place—and in our case this point $r$ will be the variable at which the last clique of size $k_1 + 1$ ends[2] The next clique after this last clique of size $k_1 + 1$ will be the clique of size $k_2 + 1$ such that only variable $r + 1$ does not belong to the previous clique of size $k_1 + 1$. The cliques in the second part of the day will now cascade as before but will be of size $k_2 + 1$.

We keep the same hyperparameters $(\alpha, \beta, \theta)$ as in Section 7.1 since they still satisfy the conditions in Section 3.1 for the given perfect order $P$ of cliques. The same $K$-fold cross-validation procedure is used to select the $(k_1, k_2, r)$ model with the lowest prediction error. We found that the $(k_1 = 14, k_2 = 4, r = 58)$ differentially banded model has the lowest prediction error when using the $IW_{P_G}$ priors. These results are consistent with those from the cross-validation in Section 7.1 and show that there is compelling evidence to suggest that there are different correlation structures during different parts of the day. For a concrete illustration of the benefits of differential banding, and for comparison purposes with [2], we also considered the task of choosing the best model for predicting the second half of the day based on the first half, as done in Section 7.1. In this case, the model with the lowest prediction error using the $IW_{P_G}$ prior Bayes estimates under the $L_1$ loss turns out to be the ($k_1 = 14$, $k_2 = 1$, $r = 55$) model. Other estimators also perform comparably but we omit the details here. The corresponding forecast error from using this model in comparison with the $k$-banded models from before is given in Figure 2. Clearly, the "differentially banded" model gives us a substantial reduction in prediction error. Using the $k$-banded model yields a 23% reduction in prediction error over the sample covariance matrix $S$ and the "differentially banded" model gives a 16% improvement over the best $k$-banded model. Our new class of models gives better prediction than the $k$-banded models for almost all time points.

---

[2]Naturally one can also extend this concept to multiple banding.

The results above highlight the performance of our estimators and their versatility in different settings. Moreover, and perhaps equally important is that taking a graphical models approach yields a much richer class of models than simple banding and achieves this flexibility in a very natural way in high-dimensional problems.

**8. Discussion.** In this paper we considered the estimation of high-dimensional covariance and precision matrices in Gaussian graphical models using a family of flexible, conjugate priors with multiple shape parameters. Existing Bayesian methods resort to either using the restrictive Diaconis–Ylvisaker conjugate prior with only one-dimensional shape parameter or using MCMC methods, both of which can be completely inadequate or at times even infeasible in very high-dimensional settings. Our objective in this paper was to overcome both of these obstacles and develop a comprehensive Bayesian solution to the problem at hand.

We derived the form of the Bayes estimators under two commonly used loss functions, adapted to graphical models, for our flexible class of priors. Another important contribution of our work is the derivation of a noninformative reference prior for $\Sigma$ and $\Omega$. Finally, we observe that our Bayes estimators have good frequentist risk properties and yield shrinkage in the eigenvalues.

The unique set of properties of the approach proposed in this paper for the estimation of large covariance matrices makes it a viable and competitive methodology. Nevertheless, there is further scope to fully assess the properties of our estimators. For instance, we would like to know more about their asymptotic properties when both $r$ and $n$ become large and we would also like to know more about the behavior of their eigenvalues. These and many other questions will be the subject of further work.

## APPENDIX

**A.1. Proof of Theorem 3.1.** To compute the mean of the $IW_{P_G}$, that is, to prove Theorem 3.1, we first need to recall the definition of the normal matrix variate and to prove the lemma below. Let us first recall the form of the matrix normal distribution as given in [31], page 79. Let $X$ be an $r \times s$ random matrix. Then we say that $X \sim N(m, \sigma_1 \otimes \sigma_2)$ where $\sigma_1$ is $r \times r$ and $\sigma_2$ is $s \times s$ if its density is of the form

$$(A.1) \quad f(x) = (2\pi)^{(rs)/2} |\sigma_1|^{-(s/2)} |\sigma_2|^{-(r/2)} \exp{-\tfrac{1}{2}} \langle \sigma_1^{-1}(x-m)\sigma_2^{-1}, (x-m) \rangle.$$

LEMMA A.1. *Let* $X \in \mathbb{R}^{r \times s}$ *follow a normal distribution* $N_{r \times s}(m, \sigma_1 \otimes \sigma_2)$ *where* $\sigma_1$ *is a positive definite* $s \times s$ *matrix and* $\sigma_2$ *is a positive definite* $r \times r$ *matrix. Let* $a = (a_{ij})_{1 \leq i,j \leq s}$ *be a given fixed* $s \times s$ *matrix. Let*

$$(\sigma_1 \otimes \sigma_2)_{ij}$$

*be the* $r \times r$ *block in the* $i$th *block row and the* $j$th *block column of the* $rs \times rs$ *covariance matrix* $(\sigma_1 \otimes \sigma_2)$ *where the rows are divided into* $s$ *sets of* $r$ *rows and so are the columns. Similarly, divide the [*$rs \times rs$ *matrix* $\mathrm{vec}\, E(X)(\mathrm{vec}\, E(X))^t]$ *matrix into* $s^2$ *blocks*

$$(\mathrm{vec}\, E(X)(\mathrm{vec}\, E(X))^t)_{ij}, \qquad 1 \leq i,j \leq s,$$

*of size* $r \times r$. *Then*

$$(A.2) \qquad E(XaX^t) = \sum_{1 \leq i,j \leq s} a_{ij}((\sigma_1 \otimes \sigma_2)_{ij} + (\mathrm{vec}(m)\,\mathrm{vec}(m)^t)_{ij})$$

$$(A.3) \qquad\qquad = mam^t + \mathrm{tr}(\sigma_1 a^t)\sigma_2.$$

PROOF. Let $X_i, i = 1, \ldots, s$ denote the columns of the $r \times s$ matrix $X$; then a straightforward calculation shows that $XaX^t = \sum_{1 \leq i,j \leq s} a_{ij} X_i X_j^t$ and therefore

$$E(XaX^t) = \sum_{1 \leq i,j \leq s} a_{ij} E(X_i X_j^t).$$

Since by [31], if $X \sim N_{r \times s}(m, \sigma_1 \otimes \sigma_2)$, then $\mathrm{vec}(X) \sim N_{rs}(\mathrm{vec}(m), \sigma_1 \otimes \sigma_2)$, it is clear that

$$E(X_i X_j^t) = \mathrm{cov}(X_i, X_j) + E(X_i)E(X_j)^t = (\sigma_1 \otimes \sigma_2)_{ij} + (\mathrm{vec}(m)\,\mathrm{vec}(m)^t)_{ij},$$

$$E(XaX^t) = \sum_{1 \leq i,j \leq s} a_{ij}((\sigma_1 \otimes \sigma_2)_{ij} + (\mathrm{vec}(m)\,\mathrm{vec}(m)^t)_{ij}),$$

which gives (A.2) and (A.3) can be verified by inspection. This completes the proof of the lemma. □

PROOF OF THEOREM 3.1. Using the distributional properties of various subblocks of an $IW_{P_G}$ matrix as given in Theorem 4.4 of [27], we will compute the expected value of the different entries of the matrix $X$. Because of the difference dependences, we shall proceed in the following order:

$$
\begin{aligned}
E(x_{\langle 2 \rangle}) & \\
E(x_{C_1 \setminus S_2, S_2}) &= E(x_{[12]} x_{\langle 2 \rangle}) \\
&= E(x_{[12]}) E(x_{\langle 2 \rangle}) = E(E(x_{[12]} | x_{[1].})) E(x_{\langle 2 \rangle}), \\
E(x_{C_1 \setminus S_2}) &= E(x_{[1].}) + E(x_{[12]} x_{\langle 2 \rangle} x_{\langle 21 \rangle}) \\
&= E(x_{[1].}) + E(x_{[12]} E(x_{\langle 2 \rangle}) x_{\langle 21 \rangle}) \\
&= E(x_{[1].}) + E(E(x_{[12]} E(x_{\langle 2 \rangle}) x_{\langle 21 \rangle} | x_{[1].})), \\
E(x_{[j]}) &= E((x_{[j]} x_{\langle j \rangle}^{-1}) x_{\langle j \rangle}) = E((x_{[j]} x_{\langle j \rangle}^{-1})) E(x_{\langle j \rangle}), \qquad j = 2, \ldots, k, \\
E(x_{[j]}) &= E(x_{[j].}) + E((x_{[j]} x_{\langle j \rangle}^{-1}) x_{\langle j \rangle} (x_{\langle j \rangle}^{-1} x_{\langle j]})) \qquad j = 2, \ldots, k, \\
&= E(x_{[j].}) + E((x_{[j]} x_{\langle j \rangle}^{-1}) E(x_{\langle j \rangle})(x_{\langle j \rangle}^{-1} x_{\langle j]})) \\
&= E(x_{[j].}) + E(E((x_{[j]} x_{\langle j \rangle}^{-1}) E(x_{\langle j \rangle})(x_{\langle j \rangle}^{-1} x_{\langle j]}) | x_{[j].})).
\end{aligned}
$$

(A.4)

Following the general formula for the expectation of an inverse Wishart distribution we have

$$
\begin{aligned}
E(x_{\langle 2 \rangle}) &= \frac{\theta_{\langle 2 \rangle}}{-(\alpha_1 + ((c_1 - s_2)/2) + \gamma_2) - ((s_2 + 1)/2)} \\
&= \frac{\theta_{\langle 2 \rangle}}{-(\alpha_1 + ((c_1 + 1)/2) + \gamma_2)}.
\end{aligned}
$$

Next,

$$
\begin{aligned}
E(x_{C_1 \setminus S_2, S_2}) &= E(x_{[12]} x_{\langle 2 \rangle}) = E(x_{[12]}) E(x_{\langle 2 \rangle}) \\
&= \theta_{[12]} \frac{\theta_{\langle 2 \rangle}}{-(\alpha_1 + ((c_1 + 1)/2) + \gamma_2)} \\
&= \frac{\theta_{C_1 \setminus S_2, S_2}}{-(\alpha_1 + ((c_1 + 1)/2) + \gamma_2)}.
\end{aligned}
$$

(A.5)

Next,

$$
\begin{aligned}
E(x_{C_1 \setminus S_2}) &= E(x_{[1].}) + E(x_{[12]} x_{\langle 2 \rangle} x_{\langle 21 \rangle}) \\
&= E(x_{[1].}) + E(x_{[12]} E(x_{\langle 2 \rangle}) x_{\langle 21 \rangle}) \\
&= E(x_{[1].}) + E(E(x_{[12]} E(x_{\langle 2 \rangle}) x_{\langle 21 \rangle} | x_{[1].})) \\
&= \frac{\theta_{[1].}}{-(\alpha_1 + ((c_1 - s_2 + 1)/2))}
\end{aligned}
$$

$$(A.6) \quad + E\Bigg(\theta_{[12]}\frac{\theta_{\langle 2\rangle}}{-(\alpha_1 + ((c_1+1)/2) + \gamma_2)}\theta_{\langle 21]}$$

$$+ \frac{1}{2}\frac{\mathrm{tr}\,\theta_{\langle 2\rangle}^{-1}\theta_{\langle 2\rangle}}{-(\alpha_1 + ((c_1+1)/2) + \gamma_2)}x_{[1]\cdot}\Bigg)$$

$$= \frac{\theta_{[1]\cdot}}{-(\alpha_1 + ((c_1 - s_2 + 1)/2))} + \frac{\theta_{C_1\setminus S_2, S_2}\theta_{\langle 2\rangle}^{-1}\theta_{S_2, C_1\setminus S_2}}{-(\alpha_1 + ((c_1 + 1)/2) + \gamma_2)}$$

$$+ \frac{s_2\theta_{[1]\cdot}}{2(\alpha_1 + ((c_1+1)/2) + \gamma_2)(\alpha_1 + ((c_1 - s_2 + 1)/2))}$$

$$= \frac{\theta_{[1]\cdot}}{-(\alpha_1 + ((c_1 - s_2 + 1)/2))}\Bigg(1 - \frac{s_2}{2(\alpha_1 + ((c_1+1)/2) + \gamma_2)}\Bigg)$$

$$+ \frac{\theta_{C_1\setminus S_2, S_2}\theta_{\langle 2\rangle}^{-1}\theta_{S_2, C_1\setminus S_2}}{-(\alpha_1 + ((c_1+1)/2) + \gamma_2)}$$

and

$$(A.7) \quad \begin{aligned} E(x_{[j\rangle}) &= E((x_{[j\rangle}x_{\langle j\rangle}^{-1})x_{\langle j\rangle}) = E((x_{[j\rangle}x_{\langle j\rangle}^{-1}))E(x_{\langle j\rangle}) \\ &= \theta_{[j\rangle}\theta_{\langle j\rangle}^{-1}E(x_{\langle j\rangle}), \end{aligned}$$

where $E(x_{\langle j\rangle})$ is given by the previous calculations (note that we are computing $E(X)$ layer by layer starting from the top).

Now, using the same type of calculations as in the computation of $E(x_{[1]\cdot})$, we have, for $j = 2, \ldots, n$,

$$\begin{aligned} E(x_{[j]}) &= E(x_{[j]\cdot}) + E(E((x_{[j\rangle}x_{\langle j\rangle}^{-1})E(x_{\langle j\rangle})(x_{\langle j\rangle}^{-1}x_{\langle j]})|x_{[j]\cdot})) \\ &= \frac{\theta_{[j]\cdot}}{-(\alpha_j + ((c_j - s_j + 1)/2))} + \theta_{[j\rangle}\theta_{\langle j\rangle}^{-1}E(x_{\langle j\rangle})\theta_{\langle j\rangle}^{-1}\theta_{\langle j]} \\ &\quad + \frac{\mathrm{tr}(\theta_{\langle j\rangle}^{-1}E(x_{\langle j\rangle}))}{-2(\alpha_j + ((c_j - s_j + 1)/2))}\theta_{[j]\cdot} \\ &= \frac{\theta_{[j]\cdot}}{-(\alpha_j + ((c_j - s_j + 1)/2))}\Big(1 + \frac{1}{2}\mathrm{tr}(\theta_{\langle j\rangle}^{-1}E(x_{\langle j\rangle}))\Big) \\ &\quad + \theta_{[j\rangle}\theta_{\langle j\rangle}^{-1}E(x_{\langle j\rangle})\theta_{\langle j\rangle}^{-1}\theta_{\langle j]}. \end{aligned}$$

(A.8)

This completes the proof. $\square$

### A.2. Proof of Theorem 4.1.

PROOF. Following (4.1), we see that the log-likelihood for $\phi$ is

$$l(\phi) = l(\sigma_{[1]\cdot}^{-1}, \sigma_{[12\rangle}, \sigma_{\langle 2\rangle}^{-1}, \sigma_{[j]\cdot}^{-1}, \sigma_{[j\rangle}\sigma_{\langle j\rangle}^{-1}, j = 2, \ldots, k)$$

$$= p \log|\sigma_{[1]\cdot}^{-1}| + p \log|\sigma_{\langle 2 \rangle}^{-1}| + \sum_{j=2}^{k} p \log|\sigma_{[j]\cdot}^{-1}| - \langle x_{[1]\cdot}, \sigma_{[1]\cdot}^{-1} \rangle - \langle x_{\langle 2 \rangle}, \sigma_{\langle 2 \rangle}^{-1} \rangle$$

(A.9)

$$- \langle \langle x_{[12\rangle} - \sigma_{[12\rangle}), \sigma_{[1]\cdot}^{-1}.(x_{[1,2\rangle} - \sigma_{[1,2\rangle})x_{\langle 2 \rangle} \rangle - \sum_{j=2}^{k} \langle x_{[j]\cdot}, \sigma_{[j]\cdot}^{-1} \rangle$$

$$- \sum_{j=2}^{k} \langle (x_{[j\rangle}x_{\langle j \rangle}^{-1} - \sigma_{[j\rangle}\sigma_{\langle j \rangle}^{-1}), \sigma_{[j]\cdot}^{-1}(x_{[j\rangle}x_{\langle j \rangle}^{-1} - \sigma_{[j\rangle}\sigma_{\langle j \rangle}^{-1})x_{\langle j \rangle} \rangle.$$

In order to obtain the information matrix $H^\phi(\phi) = E(\frac{d^2 l(\phi)}{d\phi^2})$ we differentiate the log-likelihood twice. We can then see that the Fisher's information matrix is block diagonal with blocks $H_j^\phi(\phi), j = 1, \ldots, k+1$ according to

$$\phi_1 = \sigma_{\langle 2 \rangle}^{-1}, \qquad \phi_2 = (\sigma_{[1]\cdot}^{-1}, \sigma_{[12\rangle}), \qquad \phi_{j+1} = (\sigma_{[j]\cdot}^{-1}, \sigma_{[j\rangle}\sigma_{\langle j \rangle}^{-1}), \qquad j = 2, \ldots, k.$$

In fact, since $E(x_{[12\rangle}) = \sigma_{[12\rangle}$ and $E(x_{[j\rangle}x_{\langle j \rangle}^{-1}) = \sigma_{[j\rangle}\sigma_{\langle j \rangle}^{-1}, j = 2, \ldots, k$, each $H_j^\phi(\phi), j = 2, \ldots, k+1$ is itself block diagonal and its determinant is

$$\det H^\phi(\phi) = \prod_{l=1}^{k+1} \det H_l^\phi(\phi),$$

where

(A.10)
$$\det(H_1^\phi(\phi)) = |\sigma_{\langle 2 \rangle}|^{s_2+1},$$
$$\det(H_2^\phi(\phi)) = |\sigma_{[1]\cdot}|^{-s_2}|E(X_{\langle 2 \rangle})|^{c_1-s_2}|\sigma_{[1]\cdot}|^{c_1-s_2+1},$$
$$\det(H_{j+1}^\phi(\phi)) = |\sigma_{[j]\cdot}|^{c_j-s_j+1}|\sigma_{[j]\cdot}|^{-s_j}|E(X_{\langle j \rangle})|^{c_j-s_j},$$

for $j = 2, \ldots, k$.

In the general case where $X \sim W_{Q_G}(\alpha, \beta, \sigma)$ we would have to use the expression of $E(X)$ as given in (4.21) of [27]. In the case that concerns us here, $X$ is hyper Wishart and it is well known that $E(X) \propto \sigma$ and therefore

$$E(x_{\langle j \rangle}) \propto \sigma_{\langle j \rangle}.$$

Thus, up to a constant independent of $\sigma$, we have

$$\det(H_1^\phi(\phi)) = |\sigma_{\langle 2 \rangle}|^{s_2+1},$$
$$\det(H_2^\phi(\phi)) = |\sigma_{[1]\cdot}|^{-s_2}|\sigma_{\langle 2 \rangle}|^{c_1-s_2}|\sigma_{[1]\cdot}|^{c_1-s_2+1},$$
$$\det(H_{j+1}^\phi(\phi)) = |\sigma_{[j]\cdot}|^{c_j-s_j+1}|\sigma_{[j]\cdot}|^{-s_j}|\sigma_{\langle j \rangle}|^{c_j-s_j} \qquad \text{for } j = 2, \ldots, k,$$

and therefore for $l = 1, \ldots, k+1$, $\det H_l^\phi(\phi)$ are of the form

$$\det H_l^\phi(\phi) = a_l(\phi_l)b_l(\phi_1, \ldots, \phi_{l-1})$$

where $a_l$ and $b_l$ are functions from the parameter space to $\mathbb{R}^+$. More precisely, we see that

$$a_1(\phi_1) = |\sigma_{\langle 2 \rangle}|^{s_2+1}, \qquad a_2(\phi_2) = |\sigma_{[1]}.|^{c_1+1-2s_2},$$

$$a_{l+1}(\phi_{l+1}) = |\sigma_{[l]}.|^{c_l+1-2s_l}, \qquad l = 2, \ldots, k.$$

According to the theory developed for natural exponential families by [10] and recalled in Section 2.3 of [7], the reference prior for $\phi$, with the given order of its component, is

$$(A.11) \quad \pi^\phi(\phi) = |\sigma_{\langle 2 \rangle}|^{(s_2+1)/2} |\sigma_{[1]}.|^{((c_1+1)/2)-s_2} \prod_{j=2}^{k} |\sigma_{[j]}.|^{((c_j+1)/2)-s_j}.$$

In fact, we see that the prior is independent of the order of these components and (4.3) is proved.

We now want to derive the induced prior for $\sigma$. Meticulous but relatively easy computations show that the Jacobian from $\phi$ to $\hat{\sigma}^{-1}$ is equal to

$$|\sigma_{[1]}.|^{s_2} \prod_{j=2}^{k} |\sigma_{[j]}.|^{s_j}.$$

Moreover, the Jacobian from $\hat{\sigma}^{-1}$ to $\sigma$ is known (see [32]) and equal to

$$(A.12) \qquad \prod_{j=1}^{k} (|\sigma_{[j]}.||\sigma_{\langle j \rangle}|)^{-c_j-1} \prod_{j=2}^{k} (|\sigma_{\langle j \rangle}|)^{(s_j+1)},$$

where some of the separators may be identical. Therefore the Jacobian from $\phi$ to $\sigma$ is

$$J = |\sigma_{[1]}.|^{s_2-c_1-1} \prod_{j=2}^{k} |\sigma_{[j]}.|^{s_j-c_j-1} |\sigma_{\langle 2 \rangle}|^{-c_1-1} \prod_{j=2}^{k} (|\sigma_{\langle j \rangle}|)^{(s_j-c_j)}.$$

Therefore the induced prior for $\sigma$ is

$$\pi^\sigma(\sigma) = |\sigma_{\langle 2 \rangle}|^{((s_2+1)/2)-c_1-1-c_2+s_2} |\sigma_{[1]}.|^{((c_1+1)/2)-s_2+s_2-c_1-1}$$

$$\times \prod_{j=2}^{k} |\sigma_{[j]}.|^{((c_j+1)/2)-s_j+s_j-c_j-1} \prod_{j=3}^{k} (|\sigma_{\langle j \rangle}|)^{(s_j-c_j)}$$

$$= |\sigma_{\langle 2 \rangle}|^{-((c_1+c_2)/2)-((c_1+1)/2)-((c_2+1)/2)+((s_2+1)/2)+s_2} |\sigma_{[1]}.|^{-((c_1+1)/2)}$$

$$\times \prod_{j=2}^{k} |\sigma_{[j]}.|^{-((c_j+1)/2)} \prod_{j=3}^{k} (|\sigma_{\langle j \rangle}|)^{(s_j-c_j)}$$

$$= \frac{|\sigma_{C_1}|^{-((c_1+1)/2)} \prod_{j=2}^{k} |\sigma_{C_j}|^{-((c_j+1)/2)}}{|\sigma_{S_2}|^{((c_1+c_2)/2)-s_2-((s_2+1)/2)} \prod_{j=3}^{k} |\sigma_{S_j}|^{((c_j-s_j)/2)-((s_j+1)/2)}},$$

and Theorem 4.1 is proved.   □

**Acknowledgments.**   Part of this work was done while B. Rajaratnam was a postdoctoral fellow at SAMSI and H. Massam a Research Fellow in residence at SAMSI and support is gratefully acknowledged. All three authors gratefully acknowledge support from the American Institute of Mathematics during a workshop. They would also like to thank Professor J. Huang for providing them with the Call Center data, G. Letac for the compact form (A.3) of (A.2) and D. Paul for useful discussions on model selection and cross-validation in the covariance estimation context.

B. Rajaratnam
Department of Statistics
Stanford University
Stanford, California 94305-4065
USA
E-mail: brajarat@stanford.edu

H. Massam
Department of Mathematics
  and Statistics
York University
Toronto, M3J 1P3
Canada
E-mail: massamh@yorku.ca

C. M. Carvalho
Graduate School of Business
University of Chicago
Chicago, Illinois 60637
USA
E-mail: carlos.carvalho@chicagogsb.edu